\documentclass[11pt,oneside]{article}
\usepackage{amssymb,amsmath,latexsym,amsfonts,amscd,amsthm,multirow,ctable,colortbl,mathdots}
\usepackage{fancyhdr}
\usepackage[headings]{fullpage}
\usepackage{stmaryrd}
\usepackage[all]{xy}
\usepackage{rotating}

\usepackage{setspace}

\usepackage{hyperref}

\newcommand{\bea}{\begin{eqnarray}} 
\newcommand{\eea}{\end{eqnarray}} 
\newcommand{\bee}{\begin{eqnarray*}} 
\newcommand{\eee}{\end{eqnarray*}} 
\newcommand{\al}{\begin{align*}} 
\newcommand{\eal}{\end{align*}} 
\newcommand{\be}{\begin{equation}} 
\newcommand{\ee}{\end{equation}} 
\newcommand{\eq}[1]{(\ref{#1})} 
\newcommand{\bem}{\begin{pmatrix}} 
\newcommand{\eem}{\end{pmatrix}} 

\def\a{\alpha} 
 
\def\c{\gamma}

\def\e{\epsilon}    
\def\f{\phi}

\def\im{\mathrm{Im}} 
\def\inf{\infty}

\def\m{\mu} 
\def\n{\nu}

\def\w{\omega} 
\def\p{\pi}    
\def\pa{\partial}        

\def\r{\rho}                  
\def\s{\sigma}            
\def\t{\tau} 
\def\th{\theta}

\newcolumntype{R}{ >{$}r <{$}}
\newcolumntype{C}{ >{$}c <{$}}

\def\ll{\ell}
\def\LL{\Lambda}

\newcommand{\mc}[1]{\mathcal{#1}}

\newcommand{\comment}[1]{}

\newcommand{\RR}{{\mathbb R}}
\newcommand{\CC}{{\mathbb C}}
\newcommand{\ZZ}{{\mathbb Z}}
\newcommand{\QQ}{{\mathbb Q}}
\newcommand{\HH}{{\mathbb H}}

\newcommand{\ii}{{\bf i}}

\newcommand{\tpi}{2\pi\ii}

\def\jac{\operatorname{j}}
\def\reg{\operatorname{r}}

\newcommand{\Id}{I}
\newcommand{\tr}{\operatorname{{tr}}}

\newcommand{\ex}{\operatorname{e}} 

\newcommand{\SL}{\operatorname{\textsl{SL}}}      

\newcommand{\G}{\Gamma}	
\newcommand{\g}{\gamma}	

\newcommand{\MM}{\mathbb{M}}

\newtheorem{thm}{Theorem}[section]

\newtheorem{lem}[thm]{Lemma}

\theoremstyle{definition}

\theoremstyle{remark}

\numberwithin{equation}{section}

\begin{document}

\setstretch{1.4}

\title{
    \textsc{
    	Rademacher Sums and Rademacher Series
	}
    }

\author{
	Miranda C. N. Cheng\footnote{
	Universit\'e Paris 7, UMR CNRS 7586 and LPTHE, Universit\'e Paris 6, Paris, France.
	\newline\indent\indent
	{\em E-mail:} {\tt chengm@math.jussieu.fr}
	}\\
	John F. R. Duncan\footnote{
         Department of Mathematics,
         Case Western Reserve University,
         Cleveland, OH 44106,
         U.S.A.
         \newline\indent\indent
         {\em E-mail:} {\tt john.duncan@case.edu}
                  }
}

\date{}

\maketitle

\abstract{
We exposit the construction of Rademacher sums in arbitrary weights and describe their relationship to mock modular forms. We introduce the notion of Rademacher series and describe several applications, including the determination of coefficients of Rademacher sums and a very general form of Zagier duality. We then review the application of Rademacher sums and series to moonshine both monstrous and umbral and highlight several open problems. We conclude with a discussion of the interpretation of Rademacher sums in physics.
}

\clearpage

\tableofcontents

\section{Introduction}\label{sec:intro}

Modular forms are fundamental objects in number theory which have many applications in geometry, combinatorics, string theory, and other 
branches of mathematics and physics. 
One may wonder ``what are the natural ways are to obtain a modular form?'' In general we can construct a symmetric function from a non-symmetric one by summing its images under the desired group of symmetries, although if infinite symmetry is required convergence may be a problem. A refinement of this idea, pioneered by Poincar\'e (cf. \S\ref{sec:sums:pre}), is to build in the required symmetry by summing over the images of a function $f$ that is already invariant under a (large enough) subgroup of the full group of symmetries. 
Then we may restrict the summation to representatives of cosets of the subgroup fixing $f$ and still expect to obtain a fully symmetric function. 

For instance, to obtain a modular form of even integral weight $w=2k$ we may, following Poincar\'e (cf. (\ref{eqn:sums:poi})), take $f(\t)=\ex(m\t)$ where $m$ is an integer, $\t$ is a parameter on the upper-half plane $\HH$, and here and everywhere else in the article we employ the notation
\begin{gather}
	\ex(x)=e^{2\pi\ii x}.
\end{gather}
Then the subgroup of $\G=\SL_2(\ZZ)$ leaving $f$ invariant is just the subgroup of upper-triangular matrices, which we denote $\G_{\inf}$ since its elements are precisely those that fix the infinite cusp of $\G$ (cf. \S\ref{sec:sums:pre}). Thus we are led to consider the sum
\begin{gather}\label{eqn:intro:Poi}
	\sum_{ \left(\begin{smallmatrix} a&b\\c&d\end{smallmatrix}\right) \in \G_\infty\backslash \G}   
	\ex\left(m\frac{a\t+b}{c\t+d}\right) 
	\frac{1}{(c\t+d)^{w}}, 
\end{gather}
for $w=2k$, taken over a set of representatives for the right cosets of $\G_{\inf}$ in $\G$. When $k>1$ this sum is absolutely convergent, locally uniformly for $\t\in\HH$, and thus defines a holomorphic function on $\HH$ which is invariant for the weight $w=2k$ action of $\G$ by construction. If $m\geq 0$ then it remains bounded as $\Im(\t)\to\inf$ and is thus a modular form of weight $2k$ for $\G=\SL_2(\ZZ)$. This result was obtained by Poincar\'e in \cite{Poi_FtnMdlFtnFuc}. (See \cite{Kow_PoiAncNmbThy} for a historical discussion.)

For many choices of $w$ and $m$, however (e.g. for $w\leq 2$), the infinite sum in (\ref{eqn:intro:Poi}) is not absolutely convergent (and not even conditionally convergent if $w<1$). Nontheless, we may ask if there is some way to regularise (\ref{eqn:intro:Poi}) in the case that $w\leq 2$. One solution to this problem, for the case that $w=0$, was established by Rademacher in \cite{Rad_FuncEqnModInv}. Let $J(\t)$ denote the {\em elliptic modular invariant} normalised to have vanishing constant term, so that $J(\t)$ is the unique holomorphic function on $\HH$ satisfying $J\left(\frac{a\t+b}{c\t+d}\right)=J(\t)$ whenever $\left(\begin{smallmatrix}a&b\\c&d\end{smallmatrix}\right)\in\SL_2(\ZZ)$ and also $J(\t)=q^{-1}+O(q)$ as $\Im(\t)\to \inf$ for $q=\ex(\t)$. 
\begin{gather}\label{eqn:intro:FouExpJ}
J(\t)=q^{-1} + 196884 q + 21493760 q^2 + 864299970 q^3 + \cdots
\end{gather}
In \cite{Rad_FuncEqnModInv} Rademacher established the validity of the expression
\be\label{eqn:int:Rademacher_j}
J(\t)+12=\ex(-\t) + \lim_{K\to \inf} \sum_{\substack{\left(\begin{smallmatrix} a& b \\c&d\end{smallmatrix} \right) \in \G_{\inf}\backslash\G\\0<c<K\\ -K^2 <d < K^2 }} 
\ex\left(- \frac{a\t+b}{c\t+d}\right) - \ex\left(-\frac{a}{c}\right) 
\ee
for $J(\t)$ as a conditionally convergent sum, where $\G=\SL_2(\ZZ)$, and one can recognise the right hand side of (\ref{eqn:int:Rademacher_j}) as a modification of the $w=0$ case of (\ref{eqn:intro:Poi}) with $m=-1$. This result has been generalised to other groups $\G$, and ultimately to negative (and some positive) weights, in various works, including \cite{Kno_AbIntsMdlrFns,Kno_ConstAutFrmsSuppSeries,Kno_ConstMdlrFnsI,Kno_ConstMdlrFnsII,Nie_ConstAutInts,DunFre_RSMG,Cheng2011}. (We refer to \S\ref{sec:sums} for more details.)

These regularised Poincar\'e series, which we refer to as Rademacher sums, have several important applications. Perhaps the most obvious of these is the construction of modular forms. We will see in \S\ref{sec:sums} that modular invariance sometimes but not always survives the regularisation procedure (to be described in general in \S\ref{sec:sums:reg}). More generally, a convergent Rademacher sum (cf. (\ref{eqn:sums:RSdef})) defines a mock modular form (cf. \S\ref{sec:sums:mock}); a generalisation of the notion of modular form in which the usual weight $w$ action of a discrete group $\G$ is twisted by a modular form of weight $2-w$ (cf. (\ref{eqn:sums:gtwact})).

Another application is to the computation of coefficients of modular forms. We will see in \S\ref{sec:sums}---by way of an example, cf. (\ref{eqn:sums:E2Fou})---that the Rademacher sum construction leads quite naturally to series expressions for its Fourier coefficients. This in turn leads to the notion of {\em Rademacher series}; a construction which we introduce in \S\ref{sec:series}. To a given discrete group, multiplier system and weight, the Rademacher series construction attaches, in the convergent cases, a two-dimensional grid of values. Some of these values appear as coefficients of Rademacher sums, but this typically accounts for just half of the values in the grid; the remaining values admit other interesting interpretations. For example, certain Rademacher series encode the Fourier coefficients of Eichler integrals of modular (and mock modular) forms, as we will show in \S\ref{sec:series:dualities}. The Rademacher series construction also serves to highlight a very general version of Zagier duality for Rademacher sums, whereby the set of coefficients of two families of mock modular forms in dual weights are shown to coincide, up to sign (cf. \S\ref{sec:series:dualities}).

Moreover, as we will discuss in great length in \S\ref{sec:egs}, Rademacher sums play a crucial role in the study of moonshine. 
We treat the monstrous case in \S\ref{sec:egs:mon}, the case of Mathieu moonshine in \S\ref{sec:egs:mat}, and the recently discovered umbral moonshine in \S\ref{sec:egs:umb}. We will also highlight some important open problems in this section.

Finally, an important application to physics was first pointed out in \cite{Dijkgraaf2007}. It was argued there that some Rademacher sums admit a natural physical interpretation in terms of quantum gravity 
via the so-called AdS/CFT correspondence. This interpretation has led to various work relating Rademacher sums to physical theories, and in particular to the article \cite{DunFre_RSMG} which applied the Rademacher sum construction to monstrous moonshine. One of the main results of \cite{DunFre_RSMG} is the reformulation of the genus zero property of monstrous moonshine in terms of Rademacher sums. The importance of this development has been reinforced recently by further applications \cite{Cheng2011,CheDun_M24MckAutFrms,UM}. The applications of Rademacher sums in physics will be discussed in \S\ref{sec:phys}.

\section{Rademacher Sums}\label{sec:sums}

\subsection{Preliminaries}\label{sec:sums:pre}

The group $\SL_2(\RR)$ acts naturally on the upper-half plane $\HH$ by orientation preserving isometries according to the rule 
\begin{gather}
	\bem
	a&b\\c&d
	\eem
	\t=\frac{a\t+b}{c\t+d}.
\end{gather}
For $\g\in \SL_2(\RR)$ define $\jac(\g,\t)$ to be the derivative (with respect to $\t$) of this action, so that
\begin{gather}
	\jac(\g,\t)=(c\t+d)^{-2}
\end{gather}
when $(c,d)$ is the lower row of $\g$. Let $\G$ be a subgroup of $\SL_2(\RR)$ that contains $\pm \Id$ and is commensurable with the modular group $\SL_2(\ZZ)$ and write $\G_{\infty}$ for the subgroup of $\G$ consisting of upper-triangular matrices. Then $\G_{\infty}$ is a subgroup of $\G$ isomorphic to $\ZZ\times \ZZ/2$ and is precisely the set of $\g\in\G$ for which the limit of $\g\t$ as $\Im(\t)\to\inf$ fails to be finite. (We write $\Im(\t)$ for the imaginary part of $\t$.) We set 
\begin{gather}\label{eqn:sums:T}
T=\left(\begin{matrix}1&1\\0&1\end{matrix}\right)
\end{gather}
so that $T\t=\t+1$ for $\t\in \HH$, and we write $T^h$ for $\left(\begin{smallmatrix}1&h\\0&1\end{smallmatrix}\right)$. Then there is a unique $h>0$ such that $\G_{\infty}=\langle T^h,-\Id\rangle$ and we call this $h$ the {\em width} of $\G$ at infinity. Evidently $\jac(\g,\t)=1$ for $\g\in \G_{\infty}$. 

The groups we encounter in applications typically contain and normalise the {\em Hecke congruence group} $\G_0(n)$ for some $n$.
\begin{gather}\label{eqn:sums:G0n}
	\G_0(n)=\left\{
	\begin{pmatrix}
	a&b\\c&d
	\end{pmatrix}
	\in\SL_2(\ZZ)
	\mid
	c\equiv 0\pmod{n}
	\right\}
\end{gather}
Observe that $\G_0(n)$ has width $1$ at infinity. A beautiful description of the normaliser $N(\G_0(n))$ of $\G_0(n)$ is given in \cite[\S3]{conway_norton}, and from this one can see that the width of $N(\G_0(n))$ at infinity is $1/h$ where $h$ is the largest divisor of $24$ for which $h^2$ divides $n$. 

For $w\in \RR$ say that a function $\psi:\G\to \CC$ is a {\em multiplier system} for $\G$ with weight $w$ if
\begin{gather}\label{eqn:sums:mult}
	\psi(\g_1)\psi(\g_2)\jac(\g_1,\g_2\t)^{w/2}\jac(\g_2,\t)^{w/2}
	=
	\psi(\g_1\g_2)\jac(\g_1\g_2,\t)^{w/2}
\end{gather}
for all $\g_1,\g_2\in \G$ where here and everywhere else in this paper we choose the principal branch of the logarithm (cf. (\ref{eqn:fun:pbranch})) in order to define the exponential $x\mapsto x^s$ in case $s$ is not an integer. 

Note that a multiplier system of weight $w$ is also a multiplier system of weight $w+2k$ for any integer $k$ since $\jac(\g_1,\g_2\t)\jac(\g_2,\t)=\jac(\g_1\g_2,\t)$ for any $\g_1,\g_2\in\SL_2(\RR)$. Given a multiplier system $\psi$ for $\G$ with weight $w$ we may define the {\em $(\psi,w)$-action} of $\G$ on the space $\mc{O}(\HH)$ of holomorphic functions on the upper-half plane by setting
\begin{gather}\label{eqn:sums:psiw_actn}
	(f|_{\psi,w}\g)(\t)=f(\g\t)\psi(\g)\jac(\g,\t)^{w/2}
\end{gather}
for $f\in \mc{O}(\HH)$ and $\g\in \G$. We then say that $f\in \mc{O}(\HH)$ is an {\em unrestricted modular form} with multiplier $\psi$ and weight $w$ for $\G$ in the case that $f$ is invariant for this action; i.e. $f|_{\psi,w}\g=f$ for all $\g\in\G$. Since $(-\g)\t=\g\t$ and $\jac(-\Id,\t)^{w/2}=\ex(-w/2)$ the multiplier $\psi$ must satisfy the {\em consistency condition} 
\begin{gather}\label{eqn:sums:conscond}
	\psi(-\Id)=\ex\left(\frac{w}{2}\right)
\end{gather}
in order that the corresponding space(s) of unrestricted modular forms be non-vanishing. (Recall that $\ex(x)$ is used as a shorthand for $e^{2\pi\ii x}$ throughout the article.) 

Since $\G$ is assumed to be commensurable with $\SL_2(\ZZ)$ its natural action on the boundary $\hat{\RR}=\RR\cup\{\ii\inf\}$ of $\HH$ restricts to $\hat{\QQ}=\QQ\cup\{\ii\inf\}$. The orbits of $\G$ on $\hat{\QQ}$ are called the {\em cusps} of $\G$. The quotient space
\begin{gather}\label{eqn:sums:XG}
	X_{\G}=\G\backslash\HH\cup\hat{\QQ}
\end{gather}
is naturally a compact Riemann surface (cf. e.g. \cite[\S1.5]{Shi_IntThyAutFns}). We adopt the common practice of saying that $\G$ has {\em genus zero} in case $X_{\G}$ is a genus zero surface.

We assume throughout that if $\G$ does not act transitively on $\hat{\QQ}$---i.e. if $\G$ has more than one cusp---then it is contained in a group $\tilde{\G}<\SL_2(\RR)$ that is commensurable with $\SL_2(\ZZ)$ and does act transitively on $\hat{\QQ}$, and we assume that the multiplier $\psi$ for $\G$ is of the form $\psi=\rho\tilde{\psi}$ where $\rho:\G\to\CC^{\times}$ is a morphism of groups and $\tilde{\psi}$ is a multiplier for $\tilde{\G}$. With this understanding we say that an unrestricted modular form $f$ for $\G$ with multiplier $\psi$ and weight $w$ is a {\em weak modular form} in case $f$ has at most exponential growth at the cusps of $\G$; i.e. in case there exists $C>0$ such that $(f|_{\tilde{\psi},w}\s)(\t)=O(e^{C\Im(\t)})$ as $\Im(\t)\to \inf$ for any $\s\in\tilde{\G}$. We say that $f$ is a {\em modular form} if $(f|_{\tilde{\psi},w}\s)(\t)$ remains bounded as $\Im(\t)\to \inf$ for any $\s\in\tilde{\G}$, and we say $f$ is a {\em cusp form} if $(f|_{\tilde{\psi},w}\s)(\t)\to 0$ as $\Im(\t)\to\inf$ for any $\s\in\tilde{\G}$.

If $\G$ has width $h$ at infinity then any multiplier $\psi$ for $\G$ restricts to a character on $\langle T^h\rangle<\G_{\infty}$ and so we have 
\begin{gather}\label{eqn:sums:halpha}
	\psi(T^h)=e(\a)
\end{gather}
for some $\a\in\RR$, uniquely determined subject to $0\leq \a<1$. 
Then $q^{\m}=\ex(\m\t)$ is a $\G_{\infty}$-invariant function for the $(\psi,w)$-action so long as $h\m+\a\in\ZZ$, and we may attempt to construct a $\G$-invariant function---a modular form with multiplier $\psi$ and weight $w$ for $\G$---by summing the images of $q^{\m}$ over a set of coset representatives for $\G_{\inf}$ in $\G$.
\begin{gather}\label{eqn:sums:poi}
\begin{split}
	P^{[\m]}_{\G,\psi,w}(\t)
	&=
	\sum_{\g\in\G_{\inf}\backslash\G}q^{\m}|_{\psi,w}\g\\
	&=
	\sum_{\g\in\G_{\inf}\backslash\G}
	\ex(\m\g\t)
	\psi(\g)\jac(\g,\t)^{w/2}
\end{split}
\end{gather}
This is the method that was pioneered by Poincar\'e in \cite{Poi_FtnMdlFtnFuc}. If $w>2$ then this sum (\ref{eqn:sums:poi}) converges absolutely, locally uniformly in $\t$, so that $P^{[\m]}_{\G,\psi,w}$ is a well-defined holomorphic function on $\HH$, invariant under the $(\psi,w)$-action of $\G$ by construction. Although it is not immediately obvious, $P^{[\m]}_{\G,\psi,w}$ is a weak modular form in general, a modular form in case $\m\geq 0$ and a cusp form when $\m>0$. Poincar\'e considered the special case of this construction where $\G=\SL_2(\ZZ)$, the multiplier $\psi$ is trivial and the weight $w$ is an even integer not less than $4$ in \cite{Poi_FtnMdlFtnFuc}. The more general expression (\ref{eqn:sums:poi}) was introduced by Petersson in \cite{Pet_AutFrmDtgArtPoiRhn}, and following him---Petersson called $P^{[\m]}_{\G,\psi,w}$ a ``kind of Poincar\'e series''---we call $P^{[\m]}_{\G,\psi,w}$ the {\em Poincar\'e series} of weight $w$ and index $\m$ attached to the group $\G$ and the multiplier $\psi$. 

For example, in the case that $\G$ is the modular group $\SL_2(\ZZ)$ the constant multiplier $\psi\equiv 1$ is a multiplier of weight $w=2k$ on $\G$ for any integer $k$. Taking $\m=0$ and $k>1$ we obtain the function
\begin{gather}
\begin{split}
	P^{[0]}_{\G,1,2k}(\t)
	&
	=
	\sum_{\g\in\G_{\inf}\backslash\G}
	\jac(\g,\t)^k\\
	&
	=
	1+\sum_{\substack{c,d\in\ZZ\\c>0\\(c,d)=1}}{(c\t+d)^{-2k}}
\end{split}
\end{gather}
which is the {\em Eisenstein series} of weight $2k$, often denoted $E_{2k}$, with Fourier expansion 
\begin{gather}\label{eqn:sums:Eis}
	P^{[0]}_{\G,1,2k}(\t)
	=1-\frac{4k}{B_{2k}}\sum_{n>0}\sigma_{2k-1}(n)q^n
\end{gather}
where $\sigma_p(n)$ denotes the sum of the $p$-th powers of the divisors of $n$ and $B_m$ denotes the $m$-th Bernoulli number (cf. (\ref{eqn:fun:Ber})). 
One of the main results of \cite{Pet_AutFrmDtgArtPoiRhn}---and a principal application of the Poincar\'e series construction---is that, when $w>2$, the $P^{[\m]}_{\G,\psi,w}$ for varying $\m>0$ linearly span the space of cusp forms with multiplier $\psi$ and weight $w$ for $\G$.

\subsection{Regularisation}\label{sec:sums:reg}

We may ask if there is a natural way to regularise the simple summation of (\ref{eqn:sums:poi}) in the generally divergent case when $w\leq 2$; the following method, inspired by work of {Rademacher}, is just such a procedure.  

First consider the case that $w=2$. Then the sum in (\ref{eqn:sums:poi}) is generally not absolutely convergent, but can be ordered in such a way that the result is conditionally convergent and locally uniformly so in $\t$, thus yielding a holomorphic function on $\HH$. The ordering is obtained as follows. Observe that left multiplication of a matrix $\g\in \G$ by $\pm T^h$ has no effect on the lower row of $\g$ other than to change its sign in the case of $-T^h$. So the non-trivial right-cosets of $\G_{\inf}=\langle T^h,-\Id\rangle$ in $\G$ are indexed by pairs $(c,d)$ such that $c>0$ and $(c,d)$ is the lower row of some element of $\G$. For $K>0$ we define $\G_{K,K^2}$ to be the set of elements of $\G$ having lower rows $(c,d)$ satisfying $|c|<K$ and $|d|<K^2$.
\begin{gather}
	\G_{K,K^2}=\left\{\begin{pmatrix}a&b\\c&d\end{pmatrix}\in\G\mid |c|<K,\,|d|<K^2\right\}
\end{gather}
Observe that $\G_{K,K^2}$ is a union of cosets of $\G_{\inf}$ for any $K$. Now for $\psi$ a multiplier of weight $2$ we define the {\it index} $\m$ {\em Rademacher sum} $R^{[\m]}_{\G,\psi,2}$ formally by setting
\begin{gather}
\begin{split}\label{eqn:sums:Rwt2}
	R^{[\m]}_{\G,\psi,2}(\t)
	&=
	\lim_{K\to \inf}
	\sum_{\g\in\G_{\inf}\backslash\G_{K,K^2}}
	q^{\m}|_{\psi,2}\g\\
	&=
	\lim_{K\to \inf}
	\sum_{\g\in\G_{\inf}\backslash\G_{K,K^2}}
		\ex(\m\g\t)
	\psi(\g)\jac(\g,\t),
\end{split}
\end{gather}
and we may regard $R^{[\m]}_{\G,\psi,2}(\t)$ as a holomorphic function on $\HH$ in case the limit in (\ref{eqn:sums:Rwt2}) converges locally uniformly in $\t$.

As an example we take $\G=\SL_2(\ZZ)$ and $\psi\equiv 1$ and $\m=0$ in analogy with (\ref{eqn:sums:Eis}). Then we obtain the expression
\begin{gather}
\begin{split}\label{eqn:sums:SL2wt2}
	R^{[0]}_{\G,1,2}(\t)
	&=
	\lim_{K\to \inf}
	\sum_{\g\in\G_{\inf}\backslash\G_{K,K^2}}
	\jac(\g,\t)\\
	&=
	1+\lim_{K\to \inf}\sum_{\substack{0<c<K\\-K^2<d<K^2\\(c,d)=1}}{(c\t+d)^{-2}}.
\end{split}
\end{gather}
We will show now that this expression converges. For fixed $K>0$ let $R(K)$ denote the sum in (\ref{eqn:sums:SL2wt2}) so that $R^{[0]}_{\G,1,0}=1+\lim_{K\to \inf}R(K)$. Then we have
\begin{gather}
\begin{split}\label{eqn:sums:RK1}
R(K)&=\sum_{0<c<K}c^{-2}\sum_{\substack{|d|<K^2\\(c,d)=1}}(\t+d/c)^{-2}\\
	&=\sum_{0<c<K}c^{-2}\sum_{\substack{0\leq d<c\\(c,d)=1}}\left(\sum_{|n|<K^2/c}(\t+d/c+n)^{-2}+O(c/K^2)\right)
\end{split}
\end{gather}
where the term $O(c/K^2)$ accounts for the difference between summing over $n$ such that $|d+nc|<K^2$ and summing over $n$ such that $|n|<K^2/c$, and the implied constant holds locally uniformly in $\t$. The difference between the sum over $n$ in the second line of (\ref{eqn:sums:RK1}) and its limit $\sum_{n\in\ZZ}(\t+d/c+n)^{-2}$ as $K\to \inf$ is also $O(c/K^2)$, locally uniformly for $\t\in\HH$, so we obtain
\begin{gather}
	R(K)=
	\sum_{0<c<K} 
	(-4\pi^2)c^{-2}\sum_{\substack{0\leq d<c\\(c,d)=1}}\left(\sum_{n>0}n\ex(nd/c)\ex(n\t)+O(c/K^2)\right)
\end{gather}
after an application of the Lipschitz summation formula (\ref{eqn:fun:Lipsum}) with $s=2$, $\alpha=0$. We may now estimate $\sum_{0\leq d<c}O(c/K^2)=O(c^2/K^2)$ and $\sum_{0<c<K}c^{-2}O(c^2/K^2)=O(1/K)$ and so obtain 
\begin{gather}\label{eqn:sums:RK2}
		\lim_{K\to \inf}R(K)
	=
	\lim_{K\to \inf}\sum_{0<c<K}
	(-4\pi^2)
	c^{-2}\sum_{\substack{0\leq d<c\\(c,d)=1}}\sum_{n>0}n\ex(nd/c)\ex(n\t).
\end{gather}
Now let $R'(K)$ denote the summation over $c$ in (\ref{eqn:sums:RK2}). Then $R'(K)$ is an absolutely convergent sum for fixed $K>0$ (locally uniformly so for $\t\in\HH$) and so we may reorder the terms and write
\begin{gather}\label{eqn:sums:RK3}
	R'(K)=(-4\pi^2)\sum_{n>0}n\ex(n\t)\sum_{0<c<K}c^{-2}\sum_{\substack{0\leq d<c\\(c,d)=1}}\ex\left(n\frac{d}{c}\right).
\end{gather}
The summation over $d$ in (\ref{eqn:sums:RK3}) is the sum of the $n$-th powers of the primitive $c$-th roots of unity, which is to say, it is a {\em Ramanujan sum}. The associated Dirichlet series (for fixed $n$ and varying $c$) converges absolutely for $\Re(s)>1$ and admits the explicit formula
\begin{gather}\label{eqn:sums:RSDir}
	\sum_{c>0}\sum_{\substack{0\leq d<c\\(c,d)=1}}\ex\left(n\frac{d}{c}\right)c^{-s}=n^{1-s}\frac{\s_{s-1}(n)}{\zeta(s)}
\end{gather}
in this region (cf. \cite[\S IX.1]{Siv_ClsThyArtFns}), where $\zeta(s)$ is the Riemann zeta function. Taking $s=2$ in (\ref{eqn:sums:RSDir}) we conclude that $\lim_{K\to \inf}R'(K)=\sum_{n>0}(-4\pi^2)\zeta(2)^{-1}\s_1(n)q^n$, and in particular, (\ref{eqn:sums:SL2wt2}) converges, locally uniformly for $\t\in\HH$. Applying the identity $\zeta(2)=\pi^2/6$ we obtain the Fourier expansion
\begin{gather}\label{eqn:sums:E2Fou}
	R^{[0]}_{\G,1,2}(\t)
		=1-24\sum_{n>0}\sigma_1(n)q^n
\end{gather}
and recognise $R^{[0]}_{\G,1,2}$ as the {\em quasi-modular} {Eisenstein series}, often denoted $E_2$. (Another common normalisation is $G_2=2\zeta(2)E_2$, cf. \cite[\S3.10]{Apo_MdlFnsDirSerNumThy}.)

The argument just given may be readily generalised. For example, let $\G$ be an arbitrary group commensurable with $\SL_2(\ZZ)$ that contains $-\Id$ and suppose for simplicity that $\G$ has width one at infinity. Applying a method directly similar to the above we obtain the identity $R^{[0]}_{\G,1,2}=1+\lim_{K\to \inf}R'(K)$ where now
\begin{gather}\label{eqn:sums:RK4}
	R'(K)=\sum_{n>0}(-4\pi^2)n\ex(n\t)
	\sum_{\g\in\G_{\inf}\backslash\G^{\times}_{K}/\G_{\inf}}
	\ex\left(n\frac{d}{c}\right)c^{-2}.
\end{gather}
In (\ref{eqn:sums:RK4}) we write $\G^{\times}_{K}$ for the set of elements $\g=\left(\begin{smallmatrix}a&b\\c&d\end{smallmatrix}\right)\in\G$ satisfying $0<|c|<K$
\begin{gather}\label{eqn:sums:GKdef}
	\G^{\times}_{K}=\left\{\begin{pmatrix}a&b\\c&d\end{pmatrix}\in\G\mid 0<|c|<K\right\},
\end{gather}
the summation is over a (complete and irredundant) set of representatives for the double cosets of $\G_{\inf}$ in $\G^{\times}_{K}$, and in each summand in the right most summation of (\ref{eqn:sums:RK4}) the values $c$ and $d$ are chosen so that $(c,d)$ is the lower row of the representative $\g$. Then the convergence of $R^{[0]}_{\G,1,2}$, locally uniform for $\t\in\HH$, follows in case the Dirichlet series
\begin{gather}\label{eqn:sums:Z0n}
	Z_{0,n}(s)
	=
	\lim_{K\to \inf}
	\sum_{\g\in\G_{\inf}\backslash\G^{\times}_{K}/\G_{\inf}}
	\ex\left(n\frac{d}{c}\right)c^{-2s}
\end{gather}
converges at $s=1$. This series $Z_{0,n}(s)$ is a special case of a more general construction---the {\em Kloosterman zeta function}---due to Selberg \cite{Sel_EstFouCoeffs} that we will discuss further in \S\ref{sec:series} (cf. (\ref{eqn:series:KlooZeta})). It is argued in \cite{Sel_EstFouCoeffs} that (\ref{eqn:sums:Z0n}) converges absolutely for $\Re(s)>1$; we refer to \cite{DunFre_RSMG} for a verification of the convergence of (\ref{eqn:sums:Z0n}) at $s=1$ in the case that $\G$ is commensurable with $\SL_2(\ZZ)$ and contains $-\Id$. Applying this result we obtain the convergence of $R^{[0]}_{\G,1,2}$ for such groups $\G$.

Specifying the order of summation as in (\ref{eqn:sums:Rwt2}) we may, for suitable choices of $\G$ and $\psi$, obtain conditionally convergent sums
\begin{gather}
\label{eqn:sums:Rwtwgeq1}
	R^{[\m]}_{\G,\psi,w}(\t)
	=
	\lim_{K\to \inf}
	\sum_{\g\in\G_{\inf}\backslash\G_{K,K^2}}
		\ex(\m\g\t)
	\psi(\g)\jac(\g,\t)^{w/2},
\end{gather}
converging locally uniformly for $\t\in\HH$, with weights in the range $w\geq 1$. However, the technical difficulties can be expected to increase as $w$ tends to $1$ for generally the convergence of (\ref{eqn:sums:Rwtwgeq1}) requires the convergence of a Kloosterman zeta function similar to (\ref{eqn:sums:Z0n}) at $s=w/2$, which is close to the critical line $\Re(s)=1/2$ in case $w$ is close to $1$. The convergence of some Rademacher sums with $w=3/2$ is established in \cite{Cheng2011}.
\begin{thm}[\cite{Cheng2011}]\label{thm:series:MatRadConv}
Let $\G=\G_0(n)$ for $n$ a positive integer, let $h$ be a divisor of $n$ that also divides $24$ and set $\psi=\rho_{n|h}\e^{-3}$ where $\epsilon$ and  $\rho_{n|h}$ are defined by \eq{eqn:fun:dedmlt} and (\ref{eqn:egs:mat:rhonh}). Then the Rademacher sum $R^{[1/8]}_{\G,\bar{\psi},3/2}$ converges, locally uniformly for $\t\in\HH$.
\end{thm}

In order to regularise the Poincar\'e series (\ref{eqn:sums:poi}) for weights strictly less than $1$ we require to modify the terms in the sum as well as the order in which they are taken. In general, and supposing for now that $\a\neq 0$ (cf. (\ref{eqn:sums:halpha})), we define the {\em Rademacher sum} $R^{[\m]}_{\G,\psi,w}$, for $\mu$ such that $h\mu+\a\in\ZZ$, by setting
\begin{gather}
\label{eqn:sums:RSdefanot0}
	R^{[\m]}_{\G,\psi,w}(\t)
	=
	\lim_{K\to \inf}\sum_{\g\in\G_{\inf}\backslash\G_{K,K^2}}
	\ex(\m\g\t)
	\reg_w^{[\m]}(\g,\t)\psi(\g)\jac(\g,\t)^{w/2}
\end{gather}
where $\reg^{[\m]}_w(\g,\t)$ is defined to be $1$ in case $w\geq 1$ or $\g$ is upper-triangular, and is given otherwise, in terms of the complete and lower incomplete Gamma functions (cf. (\ref{eqn:fun:gam}-\ref{eqn:fun:lowgamser})), by setting
\begin{gather}\label{eqn:sums:reg}
	\reg^{[\m]}_w(\g,\t)=\frac{1}{\G(1-w)}{\g(1-w,2\pi\ii\m(\g\t-\g\inf))}.
\end{gather}
In (\ref{eqn:sums:reg}) we write $\g\inf$ for the limit of $\g\t$ as $\t\to \ii\inf$, so $\g\inf$ is none other than $a/c$ in case $\g=\left(\begin{smallmatrix}a&b\\c&d\end{smallmatrix}\right)$ for $c\neq 0$ and is undefined when $\g\in\G_{\inf}$. We trust the reader will not be confused by the two different uses of the symbol $\g$ in (\ref{eqn:sums:reg}). Note that since we employ the principal branch of the logarithm (\ref{eqn:fun:pbranch}) everywhere in this article, and, in particular, in the definition (\ref{eqn:fun:lowgamser}) of the lower incomplete Gamma function, we should restrict $\m$ to be a non-positive real number when constructing Rademacher sums $R^{[\m]}_{\G,\psi,w}$ with $w<1$, for if $\m$ is positive then $\t\mapsto 2\pi\ii\m(\g\t-\g\inf)$ covers the left-half plane and $\reg^{[\m]}_w(\g,\t)$ can fail to be continuous with respect to $\t$. 

In the case that $w< 1$ and $\a= 0$ we need a constant term correction to the specification (\ref{eqn:sums:RSdefanot0}) so that the a complete defintion is given by 
\begin{gather}\label{eqn:sums:RSdef}
	R^{[\m]}_{\G,\psi,w}(\t)
	=
	\delta_{\a,0}\frac{1}{2}c_{\G,\psi,w}(\mu,0)
	+
	\lim_{K\to\inf}\sum_{\g\in\G_{\inf}\backslash\G_{K,K^2}}
	\ex(\m\g\t)
	\reg_w^{[\m]}(\g,\t)\psi(\g)\jac(\g,\t)^{w/2}
\end{gather}
where $c_{\G,\psi,w}(\mu,0)$ is zero in case $w\geq 1$ and is given otherwise by
\begin{gather}\label{eqn:sums:cmu0}
	c_{\G,\psi,w}(\mu,0)=\frac{1}{h}\ex\left(-\frac{w}{4}\right)\frac{(2\pi)^{2-w}(-\m)^{1-w}}{\G(2-w)}
	\lim_{K\to\inf}
	\sum_{\g\in\G_{\inf}\backslash\G^{\times}_{K}/\G_{\inf}}
	\frac{\ex(\mu\g\inf)}{c(\g)^{2-w}}
	\psi(\g)
\end{gather}
where $h$ is again the width of $\Gamma$, the lower-left-hand entry of a matrix $\g\in\SL_2(\RR)$ is denoted $c(\g)$, and 
$\G^{\times}_{K}$ is as in (\ref{eqn:sums:GKdef}). As in (\ref{eqn:sums:RK4}) the summation in (\ref{eqn:sums:cmu0}) is to be taken over a (complete and irredundant) set of representatives for the double cosets of $\G_{\inf}$ in $\G^{\times}_{K}$, chosen so that $c(\g)>0$. 
The condition $\a=0$ is necessary in order that the sum in (\ref{eqn:sums:cmu0}) not depend on the choice of representatives. As we will see in due course, the constant term correction in (\ref{eqn:sums:RSdef}) is included so as to improve the modularity of the resulting function $R^{[\m]}_{\G,\psi,w}$.

As a concrete example of a Rademacher sum with weight less than $1$ we may consider the case that $\G=\SL_2(\ZZ)$ is again the modular group, $\psi\equiv 1$ and $w=0$. Then $\g(1,x)=1-e^{-x}$ according to (\ref{eqn:fun:lowgamser}) so that when $\m=-1$ the general term in the Rademacher sum (\ref{eqn:sums:RSdef}) becomes, for $\g$ non-upper-triangular,
\begin{gather}
		\ex(\m\g\t)
	\reg_w^{[\m]}(\g,\t)\psi(\g)\jac(\g,\t)^{w/2}		
	=
		\ex(-\g\t)-\ex(-\g\inf),
\end{gather}
and we obtain
\begin{gather}\label{eqn:sums:firstradsum}
	R^{[-1]}_{\G,1,0}(\t)=	\ex(-\t)+\frac{1}{2}c_{\G,1,0}(-1,0)+\lim_{K\to\inf}\sum_{\g\in\G_{\inf}\backslash\G_{K,K^2}^{\times}}\ex(-\g\t)-\ex(-\g\inf)
\end{gather}
where the superscript $\times$ in the summation indicates a restriction to non-trivial cosets of $\G_{\inf}$. The right-hand side of (\ref{eqn:sums:firstradsum}) is in fact (but for the constant correction term) the original Rademacher sum, introduced by Rademacher in \cite{Rad_FuncEqnModInv}. Rademacher's main result in \cite{Rad_FuncEqnModInv} is that the sum 
\begin{gather}
	\ex(-\t)+\lim_{K\to\inf}\sum_{\g\in\G_{\inf}\backslash\G_{K,K^2}^{\times}}\ex(-\g\t)-\ex(-\g\inf)
\end{gather}
converges to a holomorphic function on $\HH$ that is invariant for the ($\psi\equiv 1$, $w=0$) action of the modular group and has constant term $12$ in its Fourier expansion. To calculate $c_{\G,1,0}(-1,0)$ we observe that the non-trivial double cosets of $\G_{\inf}$ in $\G=\SL_2(\ZZ)$ are represented irredundantly by the matrices $\left(\begin{smallmatrix}a&b\\c&d\end{smallmatrix}\right)$ with $c>0$ and $d$ (necessarily coprime to $c$) satisfying $0\leq d<c$. So we have
\begin{gather}\label{eqn:sums:firstradsumconst}
	c_{\G,1,0}(-1,0)=4\pi^2
		\sum_{c>0}\sum_{\substack{0\leq d<c\\(c,d)=1}}
	\ex\left(-\frac{a}{c}\right)\frac{1}{c^2}
\end{gather}
where in each term in the summation $a$ is chosen so that $ad$ is congruent to $1$ modulo $c$. Now each summation over $d$ is the sum of the primitive $c$-th roots of unity for some $c$, and so the summation over $c$ in (\ref{eqn:sums:firstradsumconst}) coincides with the special case of (\ref{eqn:sums:RSDir}) in which $n=1$ and $s=2$. So we have $c_{\G,1,0}(-1,0)=4\pi^2\zeta(2)^{-1}=24$ and thus we conclude that
\begin{gather}\label{eqn:sums:RadJ}
	R^{[-1]}_{\G,1,0}(\t)=J(\t)+24
\end{gather} 
where $J$ denotes the elliptic modular invariant (cf. (\ref{eqn:intro:FouExpJ})). We refer to \cite{Kno_RadonJPoinSerNonPosWtsEichCohom} for a nice review of Rademacher's treatment of (\ref{eqn:sums:firstradsum}). 

Generalisations of Rademacher's construction (\ref{eqn:sums:firstradsum}) have been developed by various authors, including Knopp, who attached weight $0$ Rademacher sums to various groups $\G<\SL_2(\RR)$ in \cite{Kno_ConstMdlrFnsI,Kno_ConstMdlrFnsII,Kno_AbIntsMdlrFns}, and Niebur, who established a very general convergence result for Rademacher sums of arbitrary negative weight in \cite{Nie_ConstAutInts}. 
\begin{thm}[\cite{Nie_ConstAutInts}]\label{thm:sums:NieConv}
Let $\G$ be a discrete subgroup of $\SL_2(\RR)$ having exactly one cusp. Let $\psi$ be a mulitplier for $\G$ and let $w$ be a compatible weight. If $w<0$ then the Rademacher sum $R^{[\m]}_{\G,\psi,w}$ converges for any $\m<0$ such that $h\mu+\a\in\ZZ$.
\end{thm}
We remark that the method of \cite{Nie_ConstAutInts} used to demonstrate convergence certainly applies to groups having more than one cusp.

It will develop in \S\ref{sec:series} that the convergence of Rademacher sums is generally more delicate for weights in the range $0\leq w\leq 2$ than for $|w-1|>1$. In \cite{DunFre_RSMG} it is shown that the weight $0$ Rademacher sum $R^{[\m]}_{\G,1,0}$ converges for any negative integer $\m$, for any group $\G<\SL_2(\RR)$ that is commensurable with $\SL_2(\ZZ)$ and contains $-\Id$, and certain Rademacher sums of weight $1/2$ (of relevance to Mathieu moonshine, cf. \S\ref{sec:egs:mat}) are shown to converge in \cite{Cheng2011}.
\begin{thm}[\cite{DunFre_RSMG}]\label{thm:sums:wt0RadCon}\label{thm:sums:DFConv}
Let $\G$ be a subgroup of $\SL_2(\RR)$ that is commensurable with $\SL_2(\ZZ)$ and contains $-\Id$. Then the Rademacher sum $R^{[\m]}_{\G,1,0}$ converges, locally uniformly for $\t\in\HH$, for any negative integer $\m$.
\end{thm}
\begin{thm}[\cite{Cheng2011}]\label{thm:sums:MatRadCon}\label{thm:sums:CDConv}
Let $\G=\G_0(n)$ for $n$ a positive integer, let $h$ be a divisor of $n$ that also divides $24$ and set $\psi=\rho_{n|h}\e^{-3}$
where  $\rho_{n|h}$ is defined by (\ref{eqn:egs:mat:rhonh}). Then the Rademacher sum $R^{[-1/8]}_{\G,\psi,1/2}$ converges, locally uniformly for $\t\in\HH$.
\end{thm}

\subsection{Mock Modularity}\label{sec:sums:mock}

The reader will have noticed from the examples presented so far that $\G$-invariance sometimes, but not always, survives the Rademacher regularisation procedure; the Rademacher sum $R^{[0]}_{\G,1,2}=E_2$ is not invariant when $\G=\SL_2(\ZZ)$---the Eisenstein series $E_2$ is only quasi-modular (cf. (\ref{eqn:sums:E2qmod}))---whilst the original Rademacher sum $R^{[-1]}_{\G,1,0}=J+24$ is invariant. In a word, the $\G$-invariance (with respect to the $(\psi,w)$-action) 
of a (convergent) sum $R^{[\m]}_{\G,\psi,w}$ depends upon the geometry of the group $\G$. For example, supposing that $\G$ is a subgroup of $\SL_2(\RR)$ containing $-\Id$ and commensurable with $\SL_2(\ZZ)$, the Rademacher sum $R^{[-1]}_{\G,1,0}$ fails to be $\G$-invariant exactly when $\G$ does not define a genus zero quotient of $\HH$ (i.e. when the genus of $X_{\G}$ is not zero, cf. (\ref{eqn:sums:XG})) and in this case there is a function $\omega:\G\to \CC$ such that $R^{[-1]}_{\G,1,0}(\g\t)+\omega(\g)=R^{[-1]}_{\G,1,0}(\t)$ for each $\g\in\G$ (cf. \cite[Thm. 3.4.4]{DunFre_RSMG}). The sensitivity to the genus of $\G$ in this example 
is a consequence of the choices $\psi\equiv 1$ and $w=0$, as we shall see presently. For other choices of $\psi$ and $w$ the modularity or otherwise of $R^{[\m]}_{\G,\psi,w}$ will be determined by some other geometric feature of $\G$.

In general the Rademacher regularisation defines a {\em weak mock modular form} which is a function on $\HH$ that is invariant for a certain twist of the usual $\G$-action, where the twisting is determined by a(n honest) modular form with the dual weight and inverse multiplier. More precisely, suppose that $\psi$ is a multiplier system for $\G$ with weight $w$ and $g$ is a modular form for $\G$ with the {inverse multiplier system} $\bar{\psi}:\g\mapsto \overline{\psi(\g)}$ and {\em dual} weight $2-w$. Then we can use $g$ to twist the $(\psi,w)$-action of $\G$ on $\mc{O}(\HH)$ by setting
\begin{gather}\label{eqn:sums:gtwact}
	\left(f|_{\psi,w,g}\g\right)(\t)
	=
	f(\g\t)\psi(\g)\jac(\g,\t)^{w/2}
	+(2\pi\ii)^{1-w}
	\int_{-\g^{-1}\infty}^{\ii\infty}(z+\t)^{-w}\overline{g(-\bar{z})}{\rm d}z.
\end{gather}
A {weak mock modular form} for $\G$ with multiplier $\psi$, weight $w$, and {\em shadow} $g$ is a holomorphic function $f$ on $\HH$ that is invariant for the $(\psi,w,g)$-action of $\G$ defined in (\ref{eqn:sums:gtwact}) and which has at most exponential growth at the cusps of $\G$ (i.e. there exists $C>0$ such that $(f|_{\tilde{\psi},w}\s)=O(e^{C\Im(\t)})$ for all $\s\in\tilde{\G}$ as $\Im(\t)\to\inf$ where $\tilde{\G}$ and $\tilde{\psi}$ are as in \S\ref{sec:sums:pre}). A weak mock modular form is called a {\em mock modular form} in case it is bounded at every cusp. From this point of view a (weak) modular form is a (weak) mock modular form with vanishing shadow. The notion of mock modular form developed from Zwegers' ground breaking work \cite{zwegers} on Ramanujan's mock theta functions. It is very closely related to the notion of {\em automorphic integral} which was introduced by Niebur to describe the Rademacher sums of negative weight he constructed in \cite{Nie_ConstAutInts}: an automorphic integral of weight $w$ in the sense of Niebur is a weak mock modular form whose shadow is a cusp form.

Given that convergent Rademacher sums are (weak) mock modular forms we may ask for an explicit description of the corresponding shadow functions. In fact, the Rademacher machinery itself provides such a description (cf. e.g. \cite[\S3.4]{DunFre_RSMG}, \cite[\S7]{Cheng2011}). Indeed, we can expect that the Rademacher sum $R^{[\m]}_{\G,\psi,w}$, supposing it converges, is a mock modular form whose shadow $S^{[\m]}_{\G,\psi,w}$ is also given by a Rademacher sum; namely, 
\begin{gather}\label{eqn:sums:shasum}
	S^{[\m]}_{\G,\psi,w}=\frac{(-\m)^{1-w}}{\G(1-w)}R^{[-\m]}_{\G,{\bar{\psi}},2-w}.
\end{gather}
Niebur established the identity (\ref{eqn:sums:shasum}) for arbitrary negative weights and a large class of groups.
\begin{thm}[\cite{Nie_ConstAutInts}]\label{thm:sums:NieVar}
Let $\G$ be a discrete subgroup of $\SL_2(\RR)$ having exactly one cusp. Let $\psi$ be a mulitplier for $\G$ and let $w$ be a compatible weight. If $w<0$ and $\m<0$ is such that $h\mu+\a\in\ZZ$ then the Rademacher sum $R^{[\m]}_{\G,\psi,w}$ is a weak mock modular form for $\G$ with shadow given by (\ref{eqn:sums:shasum}).
\end{thm}
Again, we remark that the method of \cite{Nie_ConstAutInts} used to demonstrate mock modularity certainly applies to groups having more than one cusp. The case that $\psi\equiv 1$ and $w=0$ in (\ref{eqn:sums:shasum}) was considered in \cite{DunFre_RSMG} and results for $w=1/2$ were established in \cite{Cheng2011}.
\begin{thm}[\cite{DunFre_RSMG}]\label{thm:sums:DFVar}
Let $\G$ be a subgroup of $\SL_2(\RR)$ that is commensurable with $\SL_2(\ZZ)$ and contains $-\Id$. Then for $\m$ a negative integer the Rademacher sum $R^{[\m]}_{\G,1,0}$ is a weak mock modular form with shadow $S^{[\m]}_{\G,1,0}$ given by (\ref{eqn:sums:shasum}). 
\end{thm}
\begin{thm}[\cite{Cheng2011}]\label{thm:sums:CDVar}
Let $\G=\G_0(n)$ for $n$ a positive integer, let $h$ be a divisor of $n$ that also divides $24$ and set $\psi=\rho_{n|h}\e^{-3}$ where $\rho_{n|h}$ is defined by (\ref{eqn:egs:mat:rhonh}). Then the Rademacher sum $R^{[-1/8]}_{\G,\psi,1/2}$ is a weak mock modular form with shadow $S^{[-1/8]}_{\G,\psi,1/2}$ given by (\ref{eqn:sums:shasum}).
\end{thm}

We can see using Theorem \ref{thm:sums:DFVar} why $R^{[-1]}_{\G,1,0}$ has to be $\G$-invariant in case $\G$ has genus zero, for the shadow $S^{[-1]}_{\G,1,0}=R^{[1]}_{\G,1,2}$ is a modular form of weight $2$ with trivial multiplier, and in fact a cusp form since it is obtained by summing images of $q=\ex(\t)$ under the weight $2$ action of $\G$. (We refer the reader to \cite{DunFre_RSMG} and \cite{Cheng2011} for more on the behavior of Rademacher sums at arbitrary cusps.) The cusp forms of weight $2$ with trivial multiplier for $\G$ are in correspondence with holomorphic $1$-forms on the Riemann surface $X_{\G}$ (cf. (\ref{eqn:sums:XG})) and the dimension of the space of holomorphic $1$-forms on a Riemann surface is equal to its genus. So if $\G$ has genus zero then $X_{\G}$ has no non-zero $1$-forms and we must have $g=S^{[-1]}_{\G\,1,0}=0$ in (\ref{eqn:sums:gtwact}).

As a second example consider the sum $R^{[0]}_{\G,1,2}$ which we found in \S\ref{sec:sums:reg} to be the Eisenstein series $E_2$ when $\G=\SL_2(\ZZ)$. To compute the right-hand side of (\ref{eqn:sums:shasum}) when $\m=0$ and $w=2$ we consider a one--parameter family of multipliers $\psi_{\delta}=\epsilon^{\delta}$, with corresponding weights $w_{\delta}=2+\delta/2$, where $\epsilon:\G\to \CC$ is the multiplier system of the Dedekind eta function (cf. (\ref{eqn:fun:Ded}-\ref{eqn:fun:eps})).
Substituting $\delta/24$ for $\m$ and $w_{\delta}=2+\delta/2$ for $w$ in (\ref{eqn:sums:shasum}) we obtain 
$-12R^{[0]}_{\G,1,0}$ in the limit as $\delta\to 0$. Recalling the definition of $\reg^{[\m]}_w(\g,\t)$ and using (\ref{eqn:sums:reg}) and (\ref{eqn:fun:lowgam}) we see that $\reg^{[0]}_{0}(\g,\t)=0$ unless $\g$ belongs to $\G_{\inf}$ in which case $\reg^{[0]}_{0}(\g,\t)=1$, so we arrive at the suggestion that the shadow of $R^{[0]}_{\G,1,2}$ should be given by $S^{[0]}_{\G,1,2}=-12R^{[0]}_{\G,1,0}\equiv -12$; that is, $R^{[0]}_{\G,1,2}$ is a mock modular form with constant shadow $-12$. Taking $g\equiv -12$ in (\ref{eqn:sums:gtwact}), and writing $R(\t)$ for $R^{[0]}_{\G,1,2}(\t)$ to ease notation, we find that
\begin{gather}
\begin{split}\label{eqn:sums:E2qmod}
	R(\t)
	=(R|_{1,2,1}\g)(\t)
	&=R(\g\t)\jac(\g,\t)
	+\frac{6\ii}{\pi}
	\int_{-\g^{-1}\inf}^{\ii\inf}(z+\t)^{-2}{\rm d}z\\
	&=R(\g\t)\jac(\g,\t)+\frac{6\ii}{\pi}\frac{1}{(\t-\g^{-1}\inf)}
\end{split}
\end{gather}
for $\g\in \G$, which is in agreement with the known quasi-modularity of $E_2$ (cf. \cite[p.69]{Apo_MdlFnsDirSerNumThy}). 

Before concluding this section we remark on an alternative approach to studying the mock modular forms we have obtained above using Rademacher sums. An equivalent and more common definition of the notion of mock modular form, more closely related to Zwegers' original treatment in \cite{zwegers},  is to say that a holomorphic function $f:\HH \to \CC$ is a weak mock modular form for the group $\G$ with multiplier $\psi$, weight $w$, and shadow $g$ if the {\em completion} of $f$, denoted $\hat f$ and defined as
\begin{gather}\label{eqn:sums:gcomp}
	\hat f(\t) = f(\t)-(2\pi\ii)^{1-w}
	\int_{-\bar \t}^{\ii\infty}(z+\t)^{-w}\overline{g(-\bar{z})}{\rm d}z,
\end{gather}
is invariant for the usual (untwisted) $(\psi,w)$-action of $\G$ (cf. (\ref{eqn:sums:psiw_actn})) on real-analytic functions on $\HH$. From (\ref{eqn:sums:gcomp}) one can check that $\hat f $ is annihilated by the differential operator
\begin{gather}
	\frac{\pa}{\pa \t} (\im \t)^w \frac{\pa}{\pa \bar \t} 
\end{gather}
and hence is a {\em harmonic weak Maa{\ss} form} of weight $w$, which is to say, $\hat{f}$ is a (non-holomorphic) modular form for $\G$ with at most exponential growth at the cusps which is also an eigenfunction for the weight $w$ Laplace operator with eigenvalue $\frac{w}{2}(1-\frac{w}{2})$. (We refer to \cite[\S5]{zagier_mock} for an exposition of this.) 
For a suitably defined Poincar\'e series (adapted to the construction of Maa{\ss} forms) the function $R^{[\m]}_{\G,\psi,w}$ may then be recovered as its {\em holomorphic part}. We refer to \cite{BringmannOno2006} for the pioneering example of this approach; further examples appear in \cite{BriOno_ArtPrpCoeHlfIntWgtMaaPoiSrs,BringmannOno2010,BriOno_CoeffHmcMaaFrms}. The harmonic weak Maa{\ss} form whose holomorphic part is $R_{\G,\e^{-3},1/2}^{[-1/8]}$ was investigated in \cite{Eguchi2009a} in the cases that $\G=\SL_2(\ZZ)$ and $\G=\G_0(2)$.

\section{Rademacher Series}\label{sec:series}

The Rademacher sums of the previous section are indexed by cosets of $\G_{\inf}$ in $\G$. In this section we consider a construction---also inspired by work of Rademacher, among others, and hinted at in the definition of the constant correction term in (\ref{eqn:sums:RSdef})---of series indexed by double coset spaces $\G_{\inf}\backslash\G^{\times}/\G_{\inf}$. 
It will develop that these series---we call them {\em Rademacher series}---
recover the Fourier coefficients of the  
Rademacher sums of the previous section, but also admit other applications, such as recovering Fourier coefficients of {false theta series} (cf. (\ref{eqn:series:EicInteta3})), and Eichler integrals of (mock) modular forms more generally (cf. (\ref{eqn:series:EicInt})). In addition, the Rademacher series construction serves to illuminate a form of {\em Zagier duality} for Rademacher sums: the coincidence (up to a root of unity depending only on $w$) of the Fourier coefficients attached to the {\em dual families}
\begin{gather}
\left\{R^{[\m]}_{\G,\psi,w}\mid h\mu+\a\in\ZZ,\,\mu<0\right\},\quad
\left\{R^{[\n]}_{\G,\bar{\psi},2-w}\mid h\n-\a\in\ZZ,\,\nu<0\right\},
\end{gather}
(cf. (\ref{eqn:sums:halpha})) of Rademacher sums.

We now detail the Rademacher series construction. Suppose as before that $\G<\SL_2(\RR)$ contains $-\Id$ and is commensurable with $\SL_2(\ZZ)$. Recall that $h>0$ is chosen so that $\G_{\inf}=\langle T^h,-\Id\rangle$ (cf. (\ref{eqn:sums:T})). Given a multiplier system $\psi$ of weight $w$ for such a group $\G$, and given also $\m,\n\in\frac{1}{h}(\ZZ-\a)$ where $\psi(T^h)=\ex(\a)$ (cf. (\ref{eqn:sums:halpha})), we define the {\em Rademacher series} $c_{\G,\psi,w}(\m,\n)$ by setting
\begin{gather}\label{eqn:series:cdef}
 	c_{\Gamma,\psi,w}(\m,\n)
	=\frac{1}{h}\lim_{K\to \inf}\sum_{\g\in\G_{\inf}\backslash\G^{\times}_K/\G_{\inf}}
	K_{\gamma,\psi}(\m,\n)B_{\gamma,w}(\m,\n)
\end{gather}
where $\G^{\times}_K$ is defined as in (\ref{eqn:sums:GKdef}) and $K_{\g,\psi}$ and $B_{\g,w}$ are given by 
\begin{gather}
	K_{\g,\psi}(\m,\n)=\label{eqn:series:Kdef}
	\ex\left( \m\frac{a}{c}\right)\ex\left( \n\frac{d}{c}\right)\psi(\g),\\
	B_{\g,w}(\m,\n)=\label{eqn:series:Bdef}
	\begin{cases}
	\ex\left(-\frac{w}{4}\right)\sum_{k\geq 0}\left(\frac{2\pi}{c}\right)^{2k+w}
	\frac{(-\m)^{k}}{k!}\frac{\n^{k+w-1}}{\G(k+w)},&w\geq 1,\\
	\ex\left(-\frac{w}{4}\right)\sum_{k\geq 0}\left(\frac{2\pi}{c}\right)^{2k+2-w}
	\frac{(-\m)^{k+1-w}}{\G(k+2-w)}\frac{\n^{k}}{k!},&w\leq 1,
	\end{cases}	
\end{gather}
in case $\g=\left(\begin{smallmatrix} a&b\\c&d\end{smallmatrix}\right)$ and $c>0$. Observe that the restriction $\m,\n\in\frac{1}{h}(\ZZ-\a)$ is necessary in order that the map $\g\mapsto K_{\g,\psi}(\m,\n)B_{\g,\psi}(\m,\n)$ descend to the double coset space $\G_{\inf}\backslash\G^{\times}_K/\G_{\inf}$; 
assuming convergence we may regard $c_{\G,\psi,w}$ as a function on the {\em grid} 
\begin{gather}\label{eqn:series:grid}
	\frac{1}{h}\ZZ\times\frac{1}{h}\ZZ-\left(\frac{\a}{h},\frac{\a}{h}\right)\subset\RR^2.
\end{gather}

Note that the convergence of the expression (\ref{eqn:series:cdef}) defining $c_{\G,\psi,w}(\m,\n)$ is not obvious when $w$ lies in the range $0\leq w\leq 2$ but is relatively easy to show for $w<0$ and $2<w$. For example, if $\G=\SL_2(\ZZ)$ and $w\geq 1$ then we have the simple estimate
\begin{gather}\label{eqn:series:csimpest}
\begin{split}
	|c_{\G,\psi,w}(\m,\n)|&\leq \sum_{\g\in\G_{\inf}\backslash\G^{\times}/\G_{\inf}}|K_{\g,\psi}(\m,\n)||B_{\g,w}(\m,\n)|\\
	&\leq \sum_{c>0}c\sum_{k\geq 0}\left(\frac{2\pi}{c}\right)^{2k+w}\frac{|\m|^k|\n|^{k+w-1}}{k!\G(k+w)}\end{split}
\end{gather}
where both $c$ and $k$ are restricted to be integers and the factor $c$ appearing between the two summations serves as a crude upper bound for the number of double cosets in $\G_{\inf}\backslash\G/\G_{\inf}$ with representatives having lower-left entry equal to $c$. Consider the result of interchanging the two summations in the right-hand side of (\ref{eqn:series:csimpest}). If $w>2$ then we obtain
\begin{gather}\label{eqn:series:csimpestswp}
\begin{split}
	\sum_{k\geq 0}(2\pi)^{2k+w}\frac{|\m|^k|\n|^{k+w-1}}{k!\G(k+w)}\sum_{c>0}c^{1-2k-w}
	&\leq
	\sum_{k\geq 0}(2\pi)^{2k+w}\frac{|\m|^k|\n|^{k+w-1}}{k!\G(k+w)}
	\frac{1}{w-2}\\
	&=
	\frac{2\pi}{w-2}
	|\m|^{(1-w)/2}
	|{\n}|^{(w-1)/2}
	I_{w-1}(4\pi|\m\n|^{1/2})
\end{split}
\end{gather}
where $I_{\a}(z)$ denotes the {\em modified Bessel function of the first kind} and we have used its series expression (\ref{eqn:fun:Bes}) in the second line of (\ref{eqn:series:csimpestswp}). 
In particular, the left-hand side of (\ref{eqn:series:csimpestswp}) is absolutely convergent for $w>2$. This verifies the coincidence of the left-hand side of (\ref{eqn:series:csimpestswp}) with the right-hand side of (\ref{eqn:series:csimpest}) and thus we obtain the absolute convergence of the Rademacher series $c_{\G,\psi,w}(\m,\n)$ for $w>2$. The case that $w<0$ is similar, and for a more general group $\G$, being a union of finitely many cosets of a finite-index subgroup of $\SL_2(\ZZ)$, the necessary adjustments to the above argument are not unduly complicated. We refer to \cite{DunFre_RSMG} for the case that $\psi$ is trivial and $w$ is an even integer. (See also Theorem \ref{thm:series:wt02conv} below.) We refer to \cite{Nie_ConstAutInts} for a treatment of the case that $w<0$.

The question of convergence is more subtle in the cases that $0\leq w\leq 2$. To establish convergence for weights in this region one has to replace the $c$ appearing between the two summations in (\ref{eqn:series:csimpest}) with a more careful estimate for the {\em Kloosterman sum}
\begin{gather}\label{eqn:series:kloosum}
	S_{\G,\psi}(\m,\n,c)=
	\sum_{\substack{\g\in\G_{\inf}\backslash\G/\G_{\inf}\\
	c(\g)=c}}K_{\g,\psi}(\m,\n).
\end{gather}
In (\ref{eqn:series:kloosum}) we again write $c(\g)$ for the lower-left entry of $\g$. A beautiful approach to analysing Kloosterman sums was pioneered by Selberg in \cite{Sel_EstFouCoeffs}. Selberg introduced the {\em Kloosterman zeta function}
\begin{gather}\label{eqn:series:KlooZeta}
	Z_{\m,\n}(s)=\sum_{\g\in\G_{\inf}\backslash\G^{\times}/\G_{\inf}}{K_{\g,\psi}(\m,\n)}{c(\g)^{-2s}}
	=\sum_{c>0}S_{\G,\psi}(\m,\n,c)c^{-2s}
\end{gather}
and demonstrated that it admits an analytic continuation that is holomorphic in the half-plane $\Re(s)>1/2$ but for finitely many poles on the real line segment $1/2<s<1$. Further, these poles are determined by the vanishing or otherwise of particular Fourier coefficients of particular cusp forms for $\G$. Using this together with the growth estimates for $Z_{\m,\n}(s)$ due to Goldfeld--Sarnak \cite{GolSar_Kloo} (see also \cite{Pri_GnlzdGolSarEst}) one may, for suitable choices of $\m$ and $\n$, obtain the convergence of the series defining $c_{\G,\psi,w}(\m,\n)$. Such an approach was first implemented by Knopp in \cite{Kno_SmlPosPowTheta,Kno_SmlPosWgt}. It was applied in \cite{DunFre_RSMG} so as to establish the convergence of $c_{\G,1,w}(\m,\n)$ in weights $w=0$ and $w=2$ for arbitrary $\G$ commensurable with $\SL_2(\ZZ)$ and arbitrary $\m,\n\in\ZZ$, and it was applied in \cite{Cheng2011} to demonstrate the convergence of $c_{\G,\psi,1/2}(\m,\n)$ for $\m=-1/8$ and $\n>0$ when $\G=\G_0(n)$ for some $n$, and $\psi$ is one of the multipliers relevant for Mathieu moonshine (cf. \S\ref{sec:egs:mat}). 

\begin{thm}[\cite{DunFre_RSMG}]\label{thm:series:wt02conv}
Let $\G$ be a subgroup of $\SL_2(\RR)$ that is commensurable with $\SL_2(\ZZ)$ and contains $-\Id$. Then the Rademacher series $c_{\G,1,0}(\m,\n)$ and $c_{\G,1,2}(\m,\n)$ converge for all $\m,\n\in\ZZ$.
\end{thm}

\begin{thm}[\cite{Cheng2011}]\label{thm:series:matconv}
Let $\G=\G_0(n)$ for $n$ a positive integer, let $h$ be a divisor of $n$ that also divides $24$ and set $\psi=\rho_{n|h}\e^{-3}$ where $\rho_{n|h}$ is defined by (\ref{eqn:egs:mat:rhonh}). Then the Rademacher series $c_{\G,\psi,1/2}(-1/8,\n)$ converges for all $\n\in\ZZ-1/8$ such that $\n>0$, and the Rademacher series $c_{\G,\bar{\psi},3/2}(1/8,\n')$ converges for all $\n'\in\ZZ+1/8$ such that $\n'>0$.
\end{thm}

At this point we may recognise the expression (\ref{eqn:sums:cmu0}), defining the constant term correction to Rademacher sums with $w=0$, as a specialisation of the Rademacher series construction (\ref{eqn:series:cdef}). In particular, we can confirm that $c_{\G,\psi,w}(\mu,0)=0$ when $w\geq 1$ and $c_{\G,\psi,w}(\mu,0)$ should not be defined unless $\a=0$. Note also that $B_{\g,w}$ can be expressed conveniently in terms of Bessel functions (cf. \S\ref{sec:fun}) in case $xy\neq 0$. For example, if $x<0$ or $y>0$ then we have
\begin{gather}\label{eqn:series:BdefI}
	B_{\g,w}(\m,\n)=
	\ex\left(-\frac{w}{4}\right)
		(-\m)^{(1-w)/2}\n^{(w-1)/2}
	\frac{2\pi}{c}
	I_{|w-1|}\left(\frac{4\pi}{c}(-\m\n)^{1/2}\right)
\end{gather}
for any weight $w\in\RR$. 
(In the case that $y<0<x$ the right-hand side of (\ref{eqn:series:BdefI}) should be multiplied by $e^{\pi\ii |w-1|}$.) 

In the remainder of this section we consider some applications of the Rademacher series. 

\subsection{Coefficients of Rademacher Sums}

Expressions like that defined by (\ref{eqn:series:cdef}--\ref{eqn:series:Bdef}) first appeared in the aforementioned work \cite{Poi_FtnMdlFtnFuc} of Poincar\'e where he considered the case that $\G=\SL_2(\ZZ)$, the multiplier $\psi$ is trivial, $w$ is an even integer greater than $2$, and $\m$ is a non-negative integer. Poincar\'e obtained an expression equivalent to $c_{\G,\psi,w}(m,n)+\delta_{m,n}$ for the Fourier coefficient of $q^n$ in $P^{[m]}_{\G,\psi,w}(\t)$, for $m$ and $n$ non-negative integers. The series of \cite{Poi_FtnMdlFtnFuc} were generalised by Petersson in \cite{Pet_AutFrmDtgArtPoiRhn}, where he obtained the analogous expression 
\begin{gather}\label{eqn:series:Poicoeff}
	P^{[\m]}_{\G,\psi,w}=q^{\m}+\sum_{\substack{h\nu+\a\in\ZZ\\\n\geq 0}}c_{\G,\psi,w}(\m,\n)q^{\n}
\end{gather}
when $\G$ is the principal congruence group $\G(N)$ (the kernel of the map $\SL_2(\ZZ)\to\SL_2(\ZZ/N\ZZ)$) for some $N$. Thus we see many instances in which the Rademacher series recover the Fourier coefficients of a Poincar\'e series.

The formula (\ref{eqn:series:Poicoeff}) was established for more general subgroups $\G<\SL_2(\ZZ)$ and for weights $w\geq 2$ in \cite{Pet_UbrEntAutFrm,Pet_UbrEntAllKlasAutFrm}, and on the strength of this, together with his result that an arbitrary modular form may be written as a linear combination of Poincar\'e series, Petersson essentially solved the problem of finding convergent series expressions for the Fourier coefficients of modular forms with weight $w\geq 2$. Using the fact that the derivative of the elliptic modular invariant $J(\t)$ is a weak modular form of weight $2$, and thus a function whose coefficients can be written in terms of the $c_{\G,\psi,w}$ according to his results, Petersson was able to derive series expressions for the coefficients of the function $J(\t)$ itself, by integration. To see such expressions consider the values $c_{\G,\psi,w}(\mu,\nu)$ for $\G=\SL_2(\ZZ)$, $\psi\equiv1 $ and $w=0$. Then $h=1$, $\alpha=0$ and $(\mu,\nu)\in\ZZ\times\ZZ$ (cf. (\ref{eqn:series:grid})). 
Observing that the non-trivial double cosets of $\G_{\inf}$ in $\G=\SL_2(\ZZ)$ are represented, irredundantly, by the matrices $\left(\begin{smallmatrix}a&b\\c&d\end{smallmatrix}\right)$ with $c>0$ and $d$ coprime to $c$ satisfying $0\leq d<c$ we find that
\begin{gather}\label{eqn:series:cforJ}
	c_{\G,1,0}\left(-1,n\right)
	=
	\sum_{\substack{c>0\\0\leq d<c\\(c,d)=1}}
	\ex\left(\frac{a+nd}{c}\right)
	n^{-{1}/{2}}
	\frac{2\pi}{c}
	I_{1}\left(\frac{4\pi}{c}n^{{1}/{2}}\right)
\end{gather}
in agreement with Petersson's formula \cite[p.202]{Pet_UbrEntAutFrm} for the $n$-th coefficient of $J(\t)$, so that, according to Rademacher's identity $R^{[-1]}_{\G,1,0}=J+24$ (cf. (\ref{eqn:sums:RadJ})), we have
\begin{gather}\label{eqn:series:RforJ}
	R^{[-1]}_{\G,1,0}(\t)=q^{-1}+\sum_{n\geq 0}c_{\G,1,0}(-1,n)q^n.
\end{gather}
In particular, the Rademacher series $c_{\G,1,0}$ recover the Fourier coefficients of the Rademacher sum $R_{\G,1,0}^{[-1]}$.

In independent work Rademacher solved the problem of providing an exact formula for the partition function \cite{Rad_PtnFn} and this furnishes another instructive example, for if $p(n)$ denotes the number of partitions of the positive integer $n$ then we have
\begin{gather}
	\frac{1}{\eta(\t)}=q^{-{1}/{24}}+\sum_{n>0}p(n)q^{n-{1}/{24}}
\end{gather}
where $\eta$ denotes the Dedekind eta function (cf. (\ref{eqn:fun:Ded})). So it suffices to compute expressions for the Fourier coefficients of the (weak) modular form $1/\eta$ of weight $-1/2$. Let $\G=\SL_2(\ZZ)$ and let $\e:\G\to \CC$ denote the multiplier system of $\eta$ (cf. (\ref{eqn:fun:eps})). Then $\bar{\e}=\e^{-1}$ is a multiplier system in weight $w=-{1}/{2}$ for $\G$ with $h=1$ and $\a={1}/{24}$  and so we may consider the values $c_{\G,\bar{\e},-1/2}(-{1}/{24},n-{1}/{24})$ for $n$ a positive integer. Comparing with the explicit formula (\ref{eqn:fun:dedmlt}) for $\e$ we find that
\begin{gather}\label{eqn:series:cforpn}
	c_{\G,\bar{\e},-{1}/{2}}\left(-\frac{1}{24},n-\frac{1}{24}\right)
	=
	\sum_{\substack{c>0\\0\leq d<c\\(c,d)=1}}
	\ex\left(n\frac{d}{c}-\frac{s(d,c)}{4}\right)
	(24n-1)^{-{3}/{4}}
	\frac{2\pi}{c}
	I_{{3}/{2}}\left(\frac{\pi}{6c}\left(24n-1\right)^{{1}/{2}}\right)
\end{gather}
which is in agreement with the formula for $p(n)$ derived in \cite{Rad_PtnFn}. (The right-hand side of (\ref{eqn:series:cforpn}) is more immediately recognised in the subsequent work \cite{RadZuc_FouCoeffMdlrFrms} which gives a general description of coefficients of modular forms of negative weight for the modular group in terms of the $c_{\G,\psi,w}$ defined above and revisits the case of $1/\eta(\t)$ as a specific example on p.455.)

Rademacher went on to determine an analogue of (\ref{eqn:series:cforpn}) for the coefficients of $J$ in \cite{Rad_FouCoeffMdlrInv}. Using a completely different method to that of \cite{Pet_UbrEntAutFrm} he independently rediscovered the formula (\ref{eqn:series:cforJ}). Rademacher's motivation for the subsequent work \cite{Rad_FuncEqnModInv}, and the introduction of the original Rademacher sum $R^{[-1]}_{\G,1,0}$ (cf. (\ref{eqn:sums:firstradsum})), was to derive the modular invariance of the function $q^{-1}+\sum_{n>0}c_{\G,1,0}(-1,n)q^n$, and thereby establish its coincidence with $J$ directly, using just the expression (\ref{eqn:series:cforJ}) for $c_{\G,1,0}(-1,n)$.

We have seen now several examples in which the series $c_{\G,\psi,w}$ serve to recover coefficients of a modular form, and a Rademacher sum in particular. In general we can expect the direct relationship 
\begin{gather}\label{eqn:series:RadFou}
	R^{[\m]}_{\G,\psi,w}(\t)=q^{\m}+\sum_{\substack{h\nu+\a\in\ZZ\\\nu\geq 0}}c_{\G,\psi,w}(\m,\nu)q^{\nu}
\end{gather}
between Rademacher sums and Rademacher series, assuming that $R^{[\m]}_{\G,\psi,w}$ and all the $c_{\G,\psi,w}(\m,\n)$ with $\n\geq 0$ are convergent. To see how this relationship can be derived we may begin by replacing $\ex(\m\g\t)$ with $\ex(\m\g\inf)\ex(\m(\g\t-\g\inf))$ in (\ref{eqn:sums:RSdef}) and rewriting $\g\t-\g\inf$ as $-c^{-1}(c\t+d)^{-1}$ in case $(c,d)$ is the lower row of $\g$. Then we may proceed in a way similar to that employed in the discussion leading to (\ref{eqn:sums:E2Fou}), applying the Lipschitz summation formula (\ref{eqn:fun:Lipsum}) (and typically also its non-absolutely convergent version, Lemma \ref{lem:fun:LipSumAnlg}) together with the fact that
\begin{gather}
\psi(\g T^h) e(\m\g T^h \inf) \jac(\g T^h,\t)^{w/2} = \psi(\g ) e(\m\g \inf) \jac(\g,\t+h)^{w/2}
\end{gather}
for $h\mu+\a\in\ZZ$, and this brings us quickly to the required expression for $R^{[\m]}_{\G,\psi,w}$ as a sum of sums over the double coset space $\G_{\inf}\backslash\G^{\times}/\G_{\inf}$. We refer to \cite{DunFre_RSMG} and \cite{Cheng2011} for detailed implementations of this approach, including careful consideration of convergence.

Since the Rademacher sum $R^{[\m]}_{\G,\psi,w}$ is precisely the Poincar\'e series $P^{[\m]}_{\G,\psi,w}$ when $w>2$ we have (\ref{eqn:series:RadFou}) for $w>2$ according to the aforementioned work of Petersson. Niebur established (\ref{eqn:series:RadFou}) for arbitrary weights $w<0$ in \cite{Nie_ConstAutInts} (and thus we have that $1/\eta$ is also a Rademacher sum---namely, $1/\eta=R^{[-1/24]}_{\G,\bar{\e},-1/2}$---according to the Rademacher's formula for $p(n)$ and the identity (\ref{eqn:series:cforpn})). We have illustrated above that the convergence of the Rademacher series $c_{\G,\psi,w}$ is more subtle in case $0\leq w\leq 2$. As we have mentioned, Petersson and Rademacher independently gave the first instance of (\ref{eqn:series:RadFou}) for $w=0$; other examples were established by Knopp in \cite{Kno_ConstMdlrFnsI,Kno_ConstMdlrFnsII,Kno_AbIntsMdlrFns}. The general case that $\G$ is commensurable with $\SL_2(\ZZ)$ and contains $-\Id$, the multiplier $\psi$ is trivial and $w=0$ was proven in \cite{DunFre_RSMG}, and examples with $w=1/2$ and $w=3/2$ were established in \cite{Cheng2011}. 
\begin{thm}[\cite{DunFre_RSMG}]\label{thm:series:wt02RadFou}
Let $\G$ be a subgroup of $\SL_2(\RR)$ that is commensurable with $\SL_2(\ZZ)$ and contains $-\Id$. Then the Fourier expansion of the Rademacher sum $R^{[\m]}_{\G,1,0}$  
is given by (\ref{eqn:series:RadFou}).
\end{thm}
\begin{thm}[\cite{Cheng2011}]\label{thm:series:MatRadFou}
Let $\G=\G_0(n)$ for $n$ a positive integer, let $h$ be a divisor of $n$ that also divides $24$ and set $\psi=\rho_{n|h}\e^{-3}$ where $\rho_{n|h}$ is defined by (\ref{eqn:egs:mat:rhonh}). Then the Fourier expansions of the Rademacher sums $R^{[-1/8]}_{\G,\psi,1/2}$ and $R^{[1/8]}_{\G,\bar{\psi},3/2}$ 
are given by (\ref{eqn:series:RadFou}).
\end{thm}
Results closely related to (\ref{eqn:series:RadFou}) for weights in the range $0<w<2$ have been established by Knopp \cite{Kno_SmlPosPowTheta,Kno_SmlPosWgt}, Pribitkin \cite{Pri_SmlPosWgt_I,Pri_SmlPosWgt_II}, and 
Bringmann--Ono \cite{BringmannOno2006,BriOno_CoeffHmcMaaFrms}.

\subsection{Dualities}\label{sec:series:dualities}

The Bessel function expression (\ref{eqn:series:BdefI}) emphasises a symmetry in $B_{\g,w}$ under the exchange of a weight $w$ with it's dual weight $2-w$; namely, $-\ex(-w/2)B_{\g,2-w}(-\n,-\m)=B_{\g,w}(\m,\n)$. Replacing $\g$ with $-\g^{-1}$ in (\ref{eqn:series:Kdef}-\ref{eqn:series:Bdef}) we observe that $\ex(w/2)K_{-\g^{-1},\bar{\psi}}(-\n,-\m)=K_{\g,w}(\m,\n)$ and $B_{-\g^{-1},w}=B_{\g,w}$, and thus we obtain the {\em Zagier duality} identity
\begin{gather}\label{eqn:series:dual}
	c_{\G,\bar{\psi},2-w}(-\nu,-\mu)
	=
	c_{\G,\psi,w}(\m,\n)
\end{gather}
in case $\mu,\nu\in\frac{1}{h}(\ZZ-\a)$ (cf. (\ref{eqn:series:grid})). This may be regarded as a generalisation of the coincidence, up to a minus sign, of coefficients in certain families of modular forms in dual weights that was observed by Zagier in \cite{Zag_TrcSngMdl}. Much of the interest in Zagier duality derives from its power to give novel interpretations to coefficients of modular forms, such as in terms of traces of singular moduli in the original example \cite{Zag_TrcSngMdl}; for other generalisations and applications we refer to \cite{Rou_ZagDtyExpBorPdtHilMdrFrms,BriOno_ArtPrpCoeHlfIntWgtMaaPoiSrs,FolOno_DltyMckThtFn,Gue_HmcMaaMdlGdsEvnIntWts,ChoCho_ZagDtyHmcMaaFmsIntWgt}.

The {duality} (\ref{eqn:series:dual}) demonstrates that dual Rademacher series---attached to mutually inverse multiplier systems in dual weights---coincide up to transposition and negation of their arguments. In other words, the vertical lines in the grid of values $(\m,\n)\mapsto c_{\G,\psi,w}(\m,\n)$ are, up to sign, the horizontal lines in the {corresponding} grid $(\m',\n')\mapsto c_{\G,\bar{\psi},2-w}(\m',\n')$ for the dual Rademacher series. 
Consequently, when considering Fourier coefficients of Rademacher sums with a given weight and multiplier system one is simultaneously considering the Fourier coefficients of Rademacher sums in the dual weight. As an application of this we see that the Rademacher series $c_{\G,\psi,w}$ encode not only the Fourier expansions of the $R^{[\m]}_{\G,\psi,w}$ but also the Fourier expansions of their shadows $S^{[\m]}_{\G,\psi,w}$. For by applying (\ref{eqn:series:RadFou}) to the formula (\ref{eqn:sums:shasum}), which relates the shadow $S^{[\m]}_{\G,\psi,w}$ of $R^{[\m]}_{\G,\psi,w}$ to the dual Rademacher sum $R^{[-\m]}_{\G,\bar{\psi},2-w}$, we obtain
\begin{gather}\label{eqn:series:ShaFou}
	S^{[\m]}_{\G,\psi,w}(\t)=\frac{(-\mu)^{1-w}}{\G(1-w)}
	\left(q^{-\m}
	\sum_{\substack{h\nu-\a\in\ZZ\\\n\geq 0}}c_{\G,\psi,w}(-\nu,\m)q^{\nu}
	\right).
\end{gather}

The {\em Eichler integral} of a cusp form $f(\t)=\sum_{\nu>0}c(\n)q^{\n}$ with weight $w$ for some group $\G$ is the function $\tilde{f}(\t)$ defined by the $q$-series 
\begin{gather}\label{eqn:series:EicInt}
\tilde{f}(\t)=\sum_{\n>0} \n^{1-w}c(\n)q^{\n}. 
\end{gather}
Let us consider the effect of transposing $\m$ with $\n$ and replacing $\g$ with $-\g^{-1}$ in (\ref{eqn:series:Kdef}-\ref{eqn:series:Bdef}). We obtain $\ex(w/2)\overline{K_{-\g^{-1},\psi}(\n,\m)}=K_{\g,\psi}(\m,\n)$ and $\ex(-w/2)\overline{B_{-\g^{-1},w}(\n,\m)\m^{1-w}}=B_{\g,w}(\m,\n)\n^{1-w}$ for $w\geq 1$, and this, together with an application of (\ref{eqn:series:dual}), leads us to the {\em Eichler duality} identity
\begin{gather}\label{eqn:series:Eicdual}
	-\overline{c_{\G,\bar{\psi},2-w}(-\m,-\n)\m^{1-w}}=c_{\G,\psi,w}(\m,\n)\n^{1-w},
\end{gather}
valid for $w\geq 1$. (A similar but slightly different expression obtains when $w< 1$.) The relation (\ref{eqn:series:Eicdual}) demonstrates another application of the Rademacher series construction: the Eichler integral of the Rademcher sum $R^{[\m]}_{\G,\psi,w}$, assuming $\a\neq 0$ or $c_{\G,\psi,w}(\m,0)=0$, is computed, up to conjugation and a scalar factor, by the Rademacher series attached to the inverse multiplier system in the dual weight. 
\begin{gather}\label{eqn:series:EicdualFou}
	\tilde{R}^{[\m]}_{\G,\psi,w}(\t)=-\overline{\m^{1-w}}\sum_{\n>0}\overline{c_{\G,\bar{\psi},2-w}(-\m,-\n)}q^{\n}
\end{gather}

As an example consider the case that $\G=\SL_2(\ZZ)$ is the modular group, $\psi=\e^{3}$ and $w=3/2$. Then $R^{[1/8]}_{\G,\e^3,3/2}$ is, up to a scalar factor, the shadow of the weak mock modular form $R^{[-1/8]}_{\G,\e^{-3},1/2}$. It is shown in \cite{Cheng2011} that $R^{[1/8]}_{\G,\e^3,3/2}=-12\eta^3$ and we have $\eta(\t)^3=\sum_{n\geq 0}(-1)^n(2n+1)q^{(2n+1)^2/8}$ according to an identity due to Euler. Thus we find that 
\begin{gather}
\tilde{R}^{[1/8]}_{\G,\e^3,3/2}(\t)=-24\sqrt{2}\sum_{n\geq 0}(-1)^nq^{(2n+1)^2/8},
\end{gather}
and applying (\ref{eqn:series:EicdualFou}) to this we obtain the beautiful formula 
\begin{gather}\label{eqn:series:EicInteta3}
	c_{\G,\e^{-3},1/2}(-1/8,-n-1/8)=
	\begin{cases}
	12(-1)^{m}&\text{if $n=\binom{m}{2}$ for some $m> 0$,}\\
	0&\text{else,}
	\end{cases}
\end{gather}
when $n\geq 0$. Compare this to the fact that the values $c_{\G,\e^{-3},1/2}(-1/8,-n-1/8)$ for $n<0$ are the coefficients of the weak mock modular form $R^{[-1/8]}_{\G,\e^{-3},1/2}$ according to Theorem \ref{thm:series:MatRadFou}. (This weak mock modular form will play a special r\^ole in \S\ref{sec:egs:mat}.) The function $\sum_{n\geq 0}(-1)^nq^{(2n+1)^2/8}$, appearing here as (a rescaling of) the Eichler integral of $\eta^3$, is one of the {\em false theta series} studied by Rogers in \cite{Rog_FalseTheta} (cf. \cite{And_RamLstNbk_III}).

\section{Moonshine}\label{sec:egs}

Some of the most fascinating and powerful applications of Rademacher sums have appeared in moonshine. 
To describe them we shall start with a short discussion of the relevant modular objects.
The study of monstrous moonshine was initiated with the realisation (cf. \cite{Tho_FinGpsModFns,Tho_NmrlgyMonsEllModFn}) that the Fourier coefficients of the elliptic modular invariant $J$ (cf. (\ref{eqn:int:Rademacher_j})) encode positive integer combinations of dimensions of irreducible representations of the monster group $\MM$.
More generally, monstrous moonshine attaches a holomorphic function $T_g=q^{-1}+\sum_{n>0}c_g(n)q^n$ on the upper-half plane to each element $g$ in the Monster group $\MM$. This association is such that the Fourier coefficients of the {\em McKay--Thompson series} $T_g$ furnish characters $g\mapsto c_g(n)$ of non-trivial representations of $\MM$ (thus the function $T_g$ depends only on the conjugacy class of $g$), and such that the $T_g$ all have the following {\em genus zero property}: 
\begin{quote}
If $\G_g$ is the invariance group of $T_g$ then the natural map $T_g:\G_{g}\backslash\HH\to \CC$ extends to an isomorphism of Riemann surfaces $X_{\G_{g}}\to \hat{\CC}$.
\end{quote}
Here $\hat{\CC}$ denotes the Riemann sphere and $X_{\G}$ is the Riemann surface $\G\backslash\HH\cup\hat{\QQ}$ (cf. (\ref{eqn:sums:XG})). We are using the weight $0$ action of $\SL_2(\RR)$ with trivial multiplier, $(f|_{1,0}\g)(\t)=f(\g\t)$ (cf. (\ref{eqn:sums:psiw_actn})), to define the invariance. 

Conway--Norton introduced the term moonshine in \cite{conway_norton} and detailed many interesting features and properties of the---at that time conjectural---correspondence $g\mapsto T_g$. An explicit monster module conjecturally realizing the $T_g$ of \cite{conway_norton} as graded-traces was constructed by Frenkel--Lepowsky--Meurman in \cite{FLMPNAS,FLMBerk,FLM}, and a beautiful proof of the Conway--Norton {\em moonshine conjectures}---that these graded traces do determine functions $T_g$ with the genus zero property formulated above---was given by Borcherds in \cite{borcherds_monstrous}. 
All that notwithstanding, a clear conceptual explanation for the genus zero property of monstrous moonshine is yet to be established. 
A step towards this goal was made in \cite{DunFre_RSMG} by employing the Rademacher sum machinery, as we shall see presently in \S\ref{sec:egs:mon}.
In particular, we will show that the genus zero property is actually equivalent to fact that $T_g$ coincides (up to a constant) with the relevant Rademacher sum (cf. \eq{eqn:egs:mon:TR}).

In \cite{Eguchi2010} a remarkable observation was made relating the elliptic genus of a $K3$ surface 
to the largest Mathieu group $M_{24}$ via a decomposition of the former into a linear combination of characters of irreducible representations of the small $N=4$ superconformal algebra. The elliptic genus is a topological invariant and for any $K3$ surface it is given by the weak Jacobi form 
\begin{gather}
	Z_{K3}(\t,z)=8
	\left(
	\left(\frac{\th_2(\t,z)}{\th_2(\t,0)}\right)^2
	+	\left(\frac{\th_3(\t,z)}{\th_3(\t,0)}\right)^2
	+	\left(\frac{\th_4(\t,z)}{\th_4(\t,0)}\right)^2
	\right)
\end{gather}
of weight $0$ and index $1$. The $\th_i$ are {Jacobi theta functions} (cf. (\ref{eqn:fun:Jactht})).
When decomposed into $N=4$ characters we obtain
\begin{align}\notag
Z_{K3}(\t,z)&=  20\, {\rm ch}^{(2)}_{\frac{1}{4},0} -2\, {\rm ch}^{(2)}_{\frac{1}{4},\frac{1}{2} } 
			+ 			\sum_{n\geq 0}
		t_{n}\, {\rm ch}^{(2)}_{\frac{1}{4}+n,\frac{1}{2}}\\ 
		& = \frac{\th_1(\t,z)^2}{\eta(\t)^3}\left(24 \,\m(\t,z) +q^{-1/8}\big(-2  + \sum_{n=1}^\inf t_n q^n \big) \right)
\end{align}
for some $t_n\in\ZZ$ where $\th_1(\t,z)$ and $\m(\t,z)$ are defined in (\ref{eqn:fun:Jactht}-\ref{eqn:fun:AppLer}).
In the above equation we write
\be\label{eqn:egs:ch}
{\rm ch}^{(\ll)}_{h,j}(\t,z) = \tr_{V^{(\ll)}_{h,j}} \left( (-1)^{J_0^3}y^{J_0^3} q^{L_0-c/24}\right)
\ee
for the Ramond sector character of the unitary highest weight representation $V^{(\ll)}_{h,j}$ of the small $N=4$ superconformal algebra with central charge $c=6(\ll-1)$. 
By inspection, the first five $t_n$ are given by $t_1=90$, $t_2=462$, $t_3=1540$, $t_4=4554$, and $t_5=11592$. The surprising observation of \cite{Eguchi2010} is that each of these $t_n$ is twice the dimension of an irreducible representation of $M_{24}$.

One is thus compelled to conjecture that every $t_{n}$ may be interpreted as the dimension of an $M_{24}$-module $K_{n-1/8}$. If we define $H(\tau)$ by requiring 
\be\label{eqn:egs:H:Jac}
 Z_{K3}(\t,z)\eta(\t)^3=\th_1(\t,z)^2(a \mu(\t,z)+H(\t))\ee
 then $a=24$ and
\be\label{eqn:egs:mat:H}
H(\t) =q^{-\frac{1}{8}}\left(-2  + \sum_{n=1}^\inf t_nq^n\right)
\ee
is a slight modification of the generating function of the $t_n$. The inclusion of the term $-2$ and the factor $q^{-1/8}$ has the effect of improving the modularity: $H(\t)$ is a weak mock modular form (cf. \S\ref{sec:sums:mock}) for $\SL_2(\ZZ)$ with multiplier $\e^{-3}$ (cf. (\ref{eqn:fun:dedmlt})), weight $1/2$, and shadow $-\frac{12}{\sqrt{2\pi}}\eta^3$ (cf. (\ref{eqn:fun:Ded})). If the $t_n$ really do encode the dimensions of $M_{24}$-modules $K_{n-1/8}$  then we can expect to obtain interesting functions $H_g(\t)$---{\em McKay--Thompson series} for $M_{24}$---by replacing $t_n$ with ${\rm tr}_{K_{n-1/8}}(g)$ in (\ref{eqn:egs:mat:H}). In other words, we should also consider
\be\label{eqn:egs:mat:Hg}
H_g(\t) =
	-2q^{-{1}/{8}}  + \sum_{n=1}^\inf {\rm tr}_{K_{n-1/8}}(g)q^{n-1/8}.
\ee

Strictly speaking, to determine $H_g$ requires knowledge of the $M_{24}$-module $K = \bigoplus_{n=1}^\inf K_{n-1/8}$ whose existence remains conjectural, but one can attempt to formulate conjectural expressions for  $H_g$ by identifying a suitably distinguishing modular property that they should satisfy. 
If the property is well-chosen then it will be strong enough for us to determine concrete expressions for the $H_g$, and compatibility between the low order terms amongst the Fourier coefficients of  $H_g$ with the character table of $M_{24}$ will serve as evidence for both the validity of $H_g$ and the existence of the module $K$. 
Exactly this was done in a series of papers, starting with \cite{Cheng2010_1}, and the independent work \cite{Gaberdiel2010}, and concluding with \cite{Gaberdiel2010a} and \cite{Eguchi2010a}. 
Despite this progress no construction of the conjectured $M_{24}$-module $K$ is yet known. To find such a construction is probably the most important open problem in Mathieu moonshine at the present time. Similar remarks also apply to the more general umbral moonshine that we will describe shortly.

The strong evidence for the conjecture that $H(\t)$ encodes the graded dimension of an $M_{24}$--module invites us to consider the $M_{24}$ analogue of the Conway--Norton moonshine conjectures---this will justify the use of the term moonshine in the $M_{24}$ setting---except that it is not immediately obvious what the analogue should be. Whilst the McKay--Thompson series $H_{g}$ is a mock modular form of weight $1/2$ on some $\G_g<\SL_2(\ZZ)$ for every $g$ in $M_{24}$ \cite{Eguchi2010a}, it is not the case that $\G_g$ is a genus zero group for every $g$, and even if it were, there is no obvious sense in which a weak mock modular form of weight $1/2$ can induce an isomorphism $X_{\G}\to \hat{\CC}$, and thus no obvious analogue of the genus zero property formulated above.
A solution to this problem---the formulation of the moonshine conjecture for $M_{24}$---was recently found in \cite{Cheng2011}.  As we shall explain in \S\ref{sec:egs:mat}, the correct analogue of the genus zero property is that the McKay--Thompson series $H_{g}$ should coincide with a certain Rademacher sum attached to its invariance group $\G_g$.

It is striking that, despite 
the very different modular properties the two sets of McKay--Thompson series $H_g$ and $T_g$  
display they can be constructed in completely analogous ways in terms of Rademacher sums. We are hence led to believe that Rademacher sums 
are an integral element of the moonshine phenomenon. 
And such a belief has in fact been instrumental in the discovery of {\em umbral moonshine} \cite{UM}, whereby a finite group $G^{(\ll)}$ and a family of vector-valued mock modular forms $H^{(\ll)}_g$ for $g\in G^{(\ll)}$ is specified for each $\ll$ in $\LL=\{2,3,4,5,7,13\}$---the set of positive integers $\ll$ such that $\ll-1$ divides $12$---and these groups $G^{(\ll)}$ and vector-valued mock modular forms $H^{(\ll)}_{g}$ are conjectured to be related in a way that we shall describe presently.

Following \cite{UM} we say that a weak Jacobi form $\f(\t,z)$ of weight $0$ and index $\ll-1$ is {\em extremal} if it admits a 
decomposition 
\begin{gather}\label{eqn:egs:umb:ext}
	\f
	=	
	a_{\frac{\ll-1}{4},0}{\rm ch}^{(\ll)}_{\frac{\ll-1}{4},0} 
	+  
	a_{\frac{\ll-1}{4},\frac{1}{2}}{\rm ch}^{(\ll)}_{\frac{\ll-1}{4},\frac{1}{2}} 
			+ \sum_{0<r<\ll}
			\sum_{\substack{n\in\ZZ\\ r^2-4\ll n< 0}}
			a_{\frac{\ll-1}{4}+n,\frac{r}{2}} {\rm ch}^{(\ll)}_{\frac{\ll-1}{4}+n,\frac{r}{2}}
\end{gather}
for some $a_{h,j}\in \CC$ where the ${\rm ch}^{(\ll)}_{h,j}$ are as in (\ref{eqn:egs:ch}). In \cite{UM} it was shown that an extremal Jacobi form is unique (up to scalar multiplication) if it exists. Moreover, it was speculated that there are no extremal Jacobi forms of index $\ll-1$ unless $\ll-1$ divides $12$, and this was shown to be true for indexes in the range $1\leq \ll-1\leq 24$.
As was discussed in detail in \cite{UM}, the above decomposition of an extremal Jacobi form $\f^{(\ll)}$  
of index $\ll-1$ leads naturally to a vector-valued mock modular form $H^{(\ll)}$ with $\ll-1$ components $H^{(\ll)}_r$, $r\in\{1,\dots,\ll-1\}$. Equivalently, the components of the vector-valued mock modular form $H^{(\ll)}=(H^{(\ll)}_r)$ are the coefficients of the theta-decomposition of the pole-free part (cf. \cite{Atish_Sameer_Don}) of a meromorphic Jacobi form of weight $1$ and index $\ll$ with a simple pole at $z=0$ that is closely related to $\phi^{(\ll)}$. 

\begin{table}[h]\begin{center} 
\caption{{ The groups of umbral moonshine.} }\label{tab:egs:umb:groups}\vspace{2mm}
\begin{tabular}{C|CCCCCCCCCCC}\toprule
\ll & 2&3&4&5&7&13\\\midrule
G^{(\ll)} & M_{24}&2.M_{12}&2.AGL_3(2) & GL_2(5)/2 &SL_2(3)& \ZZ/4\ZZ\\\bottomrule
\end{tabular} \end{center}
\end{table}

In \cite{UM} it was observed that the mock modular form $H^{(\ll)}$ obtained in this way has a close relation to a certain finite group $G^{(\ll)}$ (specified in Table \ref{tab:egs:umb:groups}) and it was conjectured that for  $\ll$ such that $\ll-1$ divides $12$ there exists a naturally defined $\ZZ\times\QQ$-graded ${G}^{(\ll)}$-module 
\begin{gather}
	K^{(\ll)}=\bigoplus_{\substack{r\in\ZZ\\0<r<\ll}}K^{(\ll)}_r=\bigoplus_{\substack{r,k\in\ZZ\\0<r<\ll}}K^{(\ll)}_{r,{k-r^2/4\ll}}
\end{gather}
such that the graded dimension of $K^{(\ll)}$ is related to the vector-valued mock modular form $H^{(\ll)}$ via
\begin{gather}\label{eqn:egs:umb:HK}
	H^{(\ll)}_{r}(\tau)=-2\delta_{r,1}q^{-1/4\ll}+\sum_{\substack{k\in\ZZ\\r^2-4k\ll<0}}\dim{K^{(\ll)}_{r,k-r^2/4\ll}}q^{k-r^2/4\ll}. 
\end{gather}
 Moreover, as in monstrous and Mathieu moonshine we expect to encounter interesting functions if we replace $\dim K^{(\ll)}_{r,k-r^2/4\ll}$ with ${\rm tr}_{K^{(\ll)}_{r,k-r^2/4\ll}}(g)$ in (\ref{eqn:egs:umb:HK}) for $g\in G^{(\ll)}$. Consider the {\em umbral McKay--Thompson series} $H^{(\ll)}_{g}=(H^{(\ll)}_{g,r})$ for $g\in G^{(\ll)}$ and $\ll\in\{2,3,4,5,7,13\}$ defined, modulo a definition of $K^{(\ll)}$, by setting
\begin{gather}\label{UM_dim_formula}
	H^{(\ll)}_{g,r}(\tau)=-2\delta_{r,1}q^{-1/4\ll}+\sum_{\substack{k\in\ZZ\\r^2-4k\ll<0}}{\rm tr}_{K^{(\ll)}_{r,k-r^2/4\ll}}(g)q^{k-r^2/4\ll}.
\end{gather}
It was conjectured in \cite{UM} that  the $G^{(\ll)}$ module $K^{(\ll)}$ has the property that all the $H^{(\ll)}_{g}$ defined above transform as vector-valued mock modular forms with specified (vector-valued) shadows. We refer to \cite{UM} for warious explicit expressions for $H^{(\ll)}_{g}$. 
The fact that all the McKay--Thompson series are mock modular forms and thus come attached with shadows is the origin of the term {\em umbral moonshine}. Notice that $G^{(2)}=M_{24}$. When $\ll=2$ the umbral moonshine conjecture stated above recovers the Mathieu moonshine conjecture relating $H(\t)$ and $M_{24}$. The Rademacher sums of relevance for umbral moonshine will be discussed  in \S\ref{sec:egs:umb}.

This series of examples clearly demonstrates the importance of Rademacher sums in understanding connections between finite groups and (mock) modular forms, and yet it seems likely that the examples presented here are not exhaustive. A complete understanding of the relationships between finite groups and mock modular forms arising from Rademacher sums would be highly desirable.

\subsection{Monstrous Moonshine}\label{sec:egs:mon}

Consider the Rademacher sums $R^{[\m]}_{\G,1,0}$ attached to groups $\G<\SL_2(\RR)$ equipped with the trivial multiplier $\psi\equiv 1$ in weight $0$, and let us specialise momentarily to the index $\mu=-1$.  
As was shown in \S\ref{sec:sums}, the formula (\ref{eqn:sums:RSdef}) for $R^{[-1]}_{\G,1,0}$ reduces in this case to 
\begin{gather}\label{eqn:egs:mon:RG}
	R^{[-1]}_{\G,1,0}(\t)
	= \ex(-\t)+\frac{1}{2} c_{\G_g,1,0}(-1,0)+
	\lim_{K\to\inf}\sum_{\G_{\inf}\backslash\G^{\times}_{<K}}
	\ex(-\g\t)-\ex(-\g\inf),
\end{gather}
As was also discussed in \S\ref{sec:sums}, it was shown in \cite{DunFre_RSMG} that the expression (\ref{eqn:egs:mon:RG}) defining $R^{[-1]}_{\G,1,0}(\t)$ converges locally uniformly in $\t$  for $\G$ commensurable with $\SL_2(\ZZ)$ and containing $-\Id$, thus yielding a holomorphic function on $\HH$. Moreover,  there is a function $\w:\G\to \CC$ such that $R^{[-1]}_{\G,1,0}(\g\t)+\w(\g)=R^{[-1]}_{\G,1,0}(\t)$ for all $\t\in\HH$, and the function $\w$ is identically zero whenever $\G$ defines a genus zero quotient of the upper-half plane.
 This last fact suggests a connection between Rademacher sums and the genus zero property of monstrous moonshine: the groups $\G_g$ are all of this specific type (commensurable with the modular group, containing $-\Id$ and having genus zero) so that $R^{[-1]}_{\G_g,1,0}$ converges and is $\G_g$-invariant for every $g\in\MM$.
Furthermore, for $\G=\G_g$ the Rademacher sum $R^{[-1]}_{\G,1,0}$ induces an isomorphism $X_{\G}\to \hat{\CC}$ (cf. (\ref{eqn:sums:XG}),  \cite{DunFre_RSMG}).
   
In fact, the connection between Rademacher sums and monstrous moonshine is even stronger.    
Given any group element $g$ of the monster, the function $T_g$ may be characterised as the unique $\G_g$-invariant holomorphic function on $\HH$ with Fourier expansion of the form $T_g(\t)=q^{-1}+O(q)$ and no poles at any non-infinite cusps of $\G_g$. In particular, the Fourier expansion (at the infinite cusp) has vanishing constant term.
It follows then that the Rademacher construction with $\mu=-1$ recovers the $T_g$ exactly, up to their constant terms, so that we have
\begin{gather}\label{eqn:egs:mon:TR}
	T_g(\t)=R^{[-1]}_{\G_g,1,0}(\t)-c_{\G_g,1,0}(-1,0)
\end{gather}
for each $g\in \MM$ according to (\ref{eqn:series:RadFou}). 
Hence we see that the Rademacher sum furnishes a uniform group-theoretic construction of the monstrous McKay--Thompson series,  
a fact that is equivalent to the genus zero property of monstrous moonshine which is yet to be fully explained. 
This leads to the expectation that a suitable physical interpretation of the Rademacher sum construction should be an integral part of a conceptual understanding of the genus zero property, and perhaps moonshine itself. 
We refer to \S\ref{sec:phys} for more on the r\^ole of Rademacher sums in physics, and to \cite[\S7]{DunFre_RSMG} for a speculative discussion of the r\^ole that physics may play in explicating monstrous moonshine.

\newpage
Given the power of Rademacher sums, one might wonder if it is possible to use them to characterise the groups $\G_g$ relevant for monstrous moonshine. At first glance this seems to be unlikely for there are many more genus zero groups\footnote{Norton, in unpublished work (cf. \cite{Cum_CngSbsGpsCmmMdlGrpGns01}), has found $616$ groups $\G$ such that $\G_{\inf}=\langle T,-\Id\rangle$, the congruence group $\G_0(N)$ is contained in $\G$ for some $N$, and the coefficients of the corresponding Rademacher sum $R^{[-1]}_{\G,1,0}$ are rational, and Cummins has shown \cite{Cum_CngSbsGpsCmmMdlGrpGns01} that $6486$ genus zero groups are obtained by dropping the condition of rationality. On the other hand, there are $194$ conjugacy classes in the monster, but the two classes of order $27$ are related by inversion and thus determine the same McKay--Thompson series. There are no other coincidences amongst the $T_g$ but there are some linear relations, and curiously, the space of functions spanned linearly by the $T_g$ for $g\in \MM$ is $163$ dimensional.} commensurable with $\SL_2(\ZZ)$ than there are functions $T_g$.
Nevertheless, a natural answer to the characterisation question is found in \cite[\S6]{DunFre_RSMG}, following earlier work \cite{ConMcKSebDiscGpsM}  by Conway--McKay--Sebbar. Following \cite{ConMcKSebDiscGpsM} we employ the natural notion of  groups of {\em $n|h$-type}, whose definition is carefully discussed in \cite[\S6]{DunFre_RSMG} and  will be suppressed here. Assuming the notion of $n|h$-type, the characterisation of \cite{DunFre_RSMG} reads as follows. A group $\G<\SL_2(\RR)$ that is of $n|h$-type and is such that $\G/\G_0(nh)$ has exponent $2$ coincides with $\G_g$ for some $g\in \MM$ if and only if 
\begin{itemize}
\item
the Rademacher sum $R^{[-1]}_{\G,1,0}$ is $\G$-invariant, and
\item
the expansion of $R^{[-1]}_{\G,1,0}$ at any cusp of $\G$ is $\G_0(nh)$-invariant.
\end{itemize}
We regard the simplicity of this formulation as further evidence that Rademacher sums have an important r\^ole to play in elucidating the nature of moonshine. (The condition that $\G/\G_0(N)$ be a group of exponent $2$ can also be formulated in terms of Rademacher sums. We refer the reader to \cite[\S6]{DunFre_RSMG} for more details.)

Finally we discuss Zagier duality for the monstrous Rademacher sums. 
So far we have only considered the Rademacher sums $R^{[\m]}_{\G,1,0}$ for $\m=-1$ but the families 
\begin{gather}\label{eqn:egs:mon:R0s}
	\left\{R^{[\m]}_{\G,1,0}\mid\mu\in\ZZ,\,\mu<0\right\}
\end{gather}
for $\G$ a monstrous group are also relevant for moonshine. Set $T^{[\m]}_{\G}=R^{[\m]}_{\G,1,0}-c_{\G,1,0}(\m,0)$ so that $T^{[-1]}_{\G}=T_g$ when $\G=\G_g$. In \cite[\S\S5,7]{DunFre_RSMG} it is argued (with detail in the case of $\G=\SL_2(\ZZ)$) that the exponential of the generating function 
$\sum_{m>0}T^{[-m]}_{\G}(\t)p^m$ furnishes the graded dimension of a certain generalised Kac--Moody algebra attached to $g$ by Carnahan in \cite{carnahan} when $\G=\G_g$ for $g\in \MM$. According to the Zagier duality (\ref{eqn:series:dual}) specialised to $w=0$ the Fourier coefficients of the family $\{T^{[-m]}_{\G}\mid m\in\ZZ,\,m>0\}$ coincide, up to a minus sign, with those of the dual family 
\begin{gather}\label{eqn:egs:mon:R2s}
	\left\{R^{[\n]}_{\G,1,2}\mid\nu\in\ZZ,\,\nu< 0\right\}. 
\end{gather}
It is interesting to observe that the subtraction of the constant terms from the $R^{[\m]}_{\G,1,0}$, which is necessary in order to obtain the functions $T^{[\m]}_{\G}$ that are of direct relevance to moonshine, has a natural reinterpretation under Zagier duality: it corresponds to the omission of the Rademacher sum $R^{[0]}_{\G,1,2}$---an Eisenstein series that fails to be modular, as was observed in \S\ref{sec:sums}---from the family (\ref{eqn:egs:mon:R2s}). 

As a final remark, we observe that the coefficients $c_{\G,1,0}$ and $c_{\G,1,2}$ are related in another way as one can see by inspecting (\ref{eqn:series:cdef}-\ref{eqn:series:Bdef}); namely, $-mc_{\G,1,2}(-m,n)=nc_{\G,1,0}(-m,n)$ for $m$ and $n$ positive integers, so the Rademacher sums $R^{[-m]}_{\G,1,2}$ dual to the functions $T^{[-m]}_{\G}=R^{[-m]}_{\G,1,0}-c_{\G,1,0}(-m,0)$ of relevance to monstrous moonshine are just their normalised derivatives,
\begin{gather}
	R^{[-m]}_{\G,1,2}=-\frac{1}{m}q\frac{{\rm d}}{{\rm d}q}T^{[-m]}_{\G}.
\end{gather}

\subsection{Mathieu Moonshine}\label{sec:egs:mat}

Consider the Rademacher sums $R^{[\m]}_{\G,\psi,w}$ with $\G=\SL_2(\ZZ)$, $\psi=\e^{-3}$ and $w=1/2$. We have $\a=1/8$ when $\psi=\e^{-3}$ so the smallest non-positive possibility for the index is $\m=-1/8$. Substituting into (\ref{eqn:sums:RSdefanot0}) we find that $R^{[-1/8]}_{\G,\e^{-3},1/2}(\t)$ is given by
\begin{gather}\label{eqn:egs:mat:RM24}
	\lim_{K\to \inf}\sum_{\substack{0<c<K\\-K^2<d<K^2\\(c,d)=1}}
		\ex\left(\frac{1}{8c(c\t+d)}+\frac{d}{8c}-\frac{3s(d,c)}{2}\right)
	\frac{-\sqrt{\ii}}{\sqrt{\pi}(c\t+d)^{w}}
	{\g\left(\frac{1}{2},\frac{-\pi\ii}{4c(c\t+d)}\right)}
\end{gather}
where $s(d,c)$ is as in (\ref{eqn:fun:dedmlt}). In deriving (\ref{eqn:egs:mat:RM24}) we have used the identities $\G(1/2)=\sqrt{\pi}$ and $\g\t-\g\inf=c^{-1}(c\t+d)^{-1}$, the latter being valid in case $(c,d)$ is the lower row of $\g\in\SL_2(\RR)$. For the Rademacher series $c_{\G,\e^{-3},1/2}$ we have
\begin{gather}\label{eqn:egs:mat:cM24}
	c_{\G,\e^{-3},\frac{1}{2}}\left(-\frac{1}{8},n-\frac{1}{8}\right)
	=
	-2\pi\sum_{\substack{c>0\\0\leq d<c\\(c,d)=1}}
	\ex\left(n\frac{d}{c}-\frac{3s(d,c)}{2}\right)
	\frac{1}{c(8n-1)^{\frac{1}{4}}}
	I_{\frac{1}{2}}\left(\frac{\pi}{2c}(8n-1)^{\frac{1}{2}}\right)
\end{gather}
according to (\ref{eqn:series:cdef}-\ref{eqn:series:BdefI}) when $n$ is a positive integer. As discussed in \S\S\ref{sec:sums},\S\ref{sec:series}, If the expressions (\ref{eqn:egs:mat:RM24}) and (\ref{eqn:egs:mat:cM24}) are convergent then the latter furnishes the Fourier expansion of the former, 
\begin{gather}\label{eqn:egs:mat:Rc}
R^{[-1/8]}_{\G,\e^{-3},1/2}(\t)=q^{-1/8}+\sum_{n>0}c_{\G,\e^{-3},1/2}\left(-1/8,n-1/8\right)q^{n-1/8}. 
\end{gather}
On the other hand, the right-hand side of (\ref{eqn:egs:mat:cM24}) appeared (up to a scalar factor) earlier in \cite{Eguchi2009a} as a proposal for an explicit formula for $t_n$. This suggests that the function $H(\t)$ may be a scalar multiple of the  Rademacher sum $R_{\G,\e^{-3},1/2}^{[-1/8]}(\t)$. In fact, more is true, for in \cite{Cheng2011} it is shown that 
for each $g\in M_{24}$ there is a character $\rho_g$ on $\G_0(n_g)$, for $n_g$ the order of $g$, such that the Rademacher sum $R^{[-1/8]}_{\G_0(n_g),\rho_g\e^{-3},1/2}$ converges, locally uniformly for $\t\in\HH$, and is related to the {McKay--Thompson series} $H_g(\t)$ by
\begin{gather}\label{eqn:egs:mat:HgRad}
	H_g(\t)=-2R_{\G_0(n_g),\rho_g\e^{-3},1/2}^{[-1/8]}(\t).
\end{gather}
The Rademacher series $c_{\G_0(n_g),\rho_g\e^{-3},1/2}$ are also shown to converge in \cite{Cheng2011}, and we recover (\ref{eqn:egs:mat:Rc}) upon taking $g$ to be the identity. In a word then, Rademacher sums furnish a uniform construction of the (candidate) $H_g$ determined earlier in \cite{Cheng2010_1,Gaberdiel2010,Gaberdiel2010a,Eguchi2010a}, which constitutes  further evidence in support of their validity. The character $\rho_g$ may be specified easily: if $n=n_g$ is the order of $g$ and $h=h_g$ is the minimal length among cycles in the cycle shape of $g$ (regarded as a permutation in the unique non-trivial permutation action on $24$ points) then $\rho_g=\rho_{n|h}$ where
\begin{gather}\label{eqn:egs:mat:rhonh}
	\rho_{n|h}
	(\g)
	=\ex\left(-\frac{cd}{nh}\right)
\end{gather}
when $(c,d)$ is the lower row of $\g\in\G_0(n)$. The fact that (\ref{eqn:egs:mat:rhonh}) defines a morphism of groups $\G_0(n_g)\to \CC^{\times}$ relies upon the result that if $h$ is a divisor of $24$ then $x^2\equiv 1 \pmod{h}$ whenever $x$ is coprime to $h$ together with the fact that all the $h_g$ for $g\in M_{24}$ are divisors of $24$. We refer to \cite{Cheng2011,CheDun_M24MckAutFrms} for more detailed discussions on the multiplier $\r_{n|h}$, as well as all the other material in this section.

As briefly mentioned before, beyond furnishing a uniform construction of the $H_g$ the result (\ref{eqn:egs:mat:HgRad}) demonstrates the correct analogue of the genus zero property that is relevant to this {\em Mathieu moonshine} relating representations of $M_{24}$ to $K3$ surfaces. The rest of this subsection will be devoted to the explanation of this fact. 
Recall that there is in this case no obvious analogue of the genus zero property which holds for the monstrous McKay--Thompson series $T_g$ since some of the groups $\G_0(n_g)$ arising in Mathieu moonshine do not define genus zero quotients of $\HH$ (viz., $n_g\in\{11,14,15,23\}$). On the other hand, from the discussion of \S\ref{sec:egs:mon} we see that the genus zero property of the $T_g$ is equivalent to the fact that they are modular functions recovered from Rademacher sums as in (\ref{eqn:egs:mon:TR}). 
Therefore, the identity (\ref{eqn:egs:mat:HgRad}) proven in \cite{Cheng2011}---the property of $H_g$ to be uniformly expressible as a Rademacher sum---serves as the natural analogue of the genus zero property that is relevant for Mathieu moonshine (modulo a proof that the $H_g$ really are the McKay--Thompson series attached to a suitably defined $M_{24}$-module $K=\bigoplus_{n>0} K_{n-1/8}$).

In more detail, we note that the identity (\ref{eqn:egs:mat:HgRad}) implies that the Rademacher sums $R_{\G,\psi,1/2}^{[-1/8]}$ with 
$\G = \G_0(n_g)$ and $\psi=\rho_{n_g|h_g}\e^{-3}$ have the special property that they are mock modular forms whose shadows lie in the one-dimensional space spanned by the cusp form $\eta^3$. This must be the case because every proposed McKay--Thompson series in Mathieu moonshine has shadow proportional to $\eta^3$, a fact that is equivalent to their relation to weak Jacobi forms generalising \eq{eqn:egs:H:Jac}.
Indeed, from (\ref{eqn:sums:shasum}) we see that the shadow of the mock modular form $ H_g(\t)$ is a weight $3/2$ modular form given by
\begin{gather}
	 -2S^{[-1/8]}_{\G_0(n_g),\rho_g\e^{-3},1/2}
	=
	-\frac{1}{\sqrt{2\pi}}R^{[1/8]}_{\G_0(n_g),{\rho_g^{-1}\e^{3}},3/2}.
\end{gather}
Moreover, it is proven in \cite{Cheng2011} that  
\begin{gather}
	 -2S^{[-1/8]}_{\G_0(n_g),\rho_g\e^{-3},1/2}
	= 
	-\frac{\chi_g}{2}\frac{1}{\sqrt{2\p}}\eta^3
\end{gather}
where $\chi_g$ denotes the number of fixed points of $g$ (in the unique non-trivial permutation representation of $M_{24}$ on $24$ points). 

As is observed in \cite{Cheng2011}, it is not typical behavior of the Rademacher sum $R_{\G,\psi,1/2}^{[-1/8]}$ to have shadow lying in this particular one-dimensional space. For $n=9$, for example---note that $9$ is not the order of an element in $M_{24}$---the shadow of the Rademacher sum $R^{[-1/8]}_{\G_0(n),\e^{-3},1/2}$ is not proportional to $\eta^3$, at least according to experimental evidence. It is natural then to ask if there is a characterisation of the modular groups and the multipliers of the McKay--Thompson series $H_g$ expressible in terms of Rademacher sums, in analogy with that of \cite{DunFre_RSMG} (derived following \cite{ConMcKSebDiscGpsM}) for the monstrous case as discussed in \S\ref{sec:egs:mon}. In such a characterisation the pairs $(\G_0(n),\rho_{n|h})$ for $h$ a divisor of $n$ dividing $24$ would replace the groups of $n|h$-type, and the condition 
\begin{itemize}
\item
{ the Rademacher sum $-2R^{[-1/8]}_{\G_0(n),\rho_{n|h}\e^{-3},1/2}$ has shadow proportional to $\eta^3$}
\end{itemize}
would replace the $\G$-invariance condition in \S\ref{sec:egs:mon}. So far we do not know of any examples that do not arise as $H_g$ for some $g\in M_{24}$. It would be very interesting to determine whether or not the above conditions are sufficient to characterise the McKay--Thompson series of Mathieu moonshine.

\subsection{Umbral Moonshine}\label{sec:egs:umb}

In \S\S\ref{sec:sums},\ref{sec:series} we have described a  regularisation procedure attaching Rademacher sums $R^{[\m]}_{\G,\psi,w}$ to a group $\G<\SL_2(\RR)$, a multiplier $\psi$ for $\G$, a compatible weight $w$ and a compatible index $\m$.  
This procedure can be generalised to the vector-valued case with a higher-dimensional $\psi$ and $\m$. 
To be precise, we suppose that $\psi=(\psi_{ij})$ is a matrix-valued multiplier system, satisfying (\ref{eqn:sums:mult}) as before, for some weight $w$, and we suppose also that $\psi_{ij}(T^h)=\delta_{ij}\ex(\a_i)$ for some $0<\a_i<1$ where $h$ is such that $\G_{\inf}=\langle T^h,-\Id\rangle$. Then to a vector-valued index $\mu=[\m_i]$ such that $h\mu_i+\a_i\in\ZZ$ for all $i$ (and $\mu_i<0$ in case $w<1$) we attach the {\em (row) vector-valued Rademacher sum} 
\begin{gather}\label{eqn:egs:umb:vvR}
	R^{[\m]}_{\G,\psi,w}(\t)=\lim_{K\to\inf}\sum_{\G_{\inf}\backslash\G_{K,K^2}}\ex(\mu\g\t)\reg^{[\m]}_w(\g,\t)\psi(\g)\jac(\g,\t)^{w/2}
\end{gather}
where $\ex(\m\g\t)$ now denotes the (row) vector-valued function whose $i$-th component is $\ex(\m_i\g\t)$ and $\reg^{[\m]}_w(\g,\t)$ denotes the diagonal matrix-valued function whose $(i,i)$-th entry is $\reg^{[\m_i]}_{w}(\g,\t)$ (cf. (\ref{eqn:sums:reg})). For the sake of simplicity we exclude the case that some $\a_i=0$ in (\ref{eqn:egs:umb:vvR}). In such a case one can expect constant term corrections analogous to (\ref{eqn:sums:RSdef}).

In order to apply the above construction to the vector-valued mock modular forms relevant for umbral moonshine we have to specify the appropriate (matrix-valued) multiplier. Recall that the vector-valued mock modular forms $H^{(\ll)}$ are obtained from the decomposition of extremal Jacobi forms into $N=4$ characters. As is explained in detail in \cite{UM}, the relation to the weak Jacobi form immediately implies that the mock modular form $H^{(\ll)}$ has shadow (proportional to) $S^{(\ll)}=(S^{(\ll)}_r)$, whose components are the {\em unary theta series}
\begin{gather}\label{eqn:egs:umb:S}
	S^{(\ll)}_r(\t)=\sum_{k\in\ZZ}(2\ll k+r)q^{\frac{(2\ll k+r)^2}{4\ll}}, 
\end{gather}
while the extremality condition implies that $H^{(\ll)}$ has a single polar (non-vanishing as $\t \to i\inf$) term $-2q^{-\frac{1}{4\ll}}$ in its first component $H^{(\ll)}_1$ (cf. \eq{eqn:egs:umb:HK}). Notice that in the case that $\ll=2$ we have $S^{(2)}=(S^{(2)}_1)=(\eta^3)$ by an identity due to Euler, and this is in part a reflection of the fact that $S^{(\ll)}$ is a (vector-valued) cusp form of weight $3/2$ for $\SL_2(\ZZ)$ for all $\ll\geq 2$.

 Let $\s^{(\ll)}=(\s^{(\ll)}_{ij})$ be the multiplier system for $S^{(\ll)}$. Then from the above discussion, we would like to consider the $(\ll-1)$-vector-valued Rademacher sum $R^{[\m]}_{\G,\psi^{(\ll)},1/2}$ where $\G=\SL_2(\ZZ)$, we take $\psi^{(\ll)}$ to be the inverse of $\s^{(\ll)}$, and where we set $\mu=\mu^{(\ll)}=(-\frac{1}{4\ll},0,\ldots,0)$. As was uncovered in \cite{UM}, the Rademacher sum $R^{[\m]}_{\G,\psi^{(\ll)},1/2}$ (denoted $R^{(\ll)}$ in \cite{UM}) has special properties when $\ll-1$ is a divisor of $12$. First, in these cases $R^{[\m]}_{\G,\psi^{(\ll)},1/2}$ turns out to be a vector-valued mock modular form with shadow proportional to the vector-valued cusp form $S^{(\ll)}$ defined in (\ref{eqn:egs:umb:S}). This means that, for a suitably chosen constant $C^{(\ll)}$, the vector-valued function $R^{[\m]}_{\G,\psi^{(\ll)},1/2}$ is invariant for the $(\psi,w,G)$-action of $\G=\SL_2(\ZZ)$ on $(\ll-1)$-vector-valued holomorphic functions $F(\t)=(F_1(\t),\ldots,F_{\ll-1}(\t))$ defined, in direct analogy with (\ref{eqn:sums:gtwact}), by setting
\begin{gather}\label{eqn:egs:umb:gtwact}
	\left(F|_{\psi,w,G}\g\right)(\t)
	=
		F(\g\t)\psi(\g)\jac(\g,\t)^{w/2}
	+(2\pi\ii)^{1-w}
		\int_{-\g^{-1}\infty}^{\ii\infty}(z+\t)^{-w}\overline{G(-\bar{z})}{\rm d}z,
\end{gather}
when $\psi=\psi^{(\ll)}$, $w=1/2$, and $G(\t)=C^{(\ll)}S^{(\ll)}(\t)=C^{(\ll)}(S^{(\ll)}_1(\t),\ldots,S^{(\ll)}_{\ll-1}(\t))$. Second, it appears to have a close relation to the group $G^{(\ll)}$ as described in \eq{UM_dim_formula}.

As the reader might have noticed, in case $\ll=2$ the function $R^{[\m]}_{\G,\psi^{(\ll)},1/2}$ has a single component which by definition coincides with $R^{[-1/8]}_{\G,\e^{-3},1/2}$. Thus $-2R^{[\m]}_{\G,\psi^{(\ll)},1/2}$ recovers the mock modular form $H(\t)$ of importance in Mathieu moonshine (and discussed in \S\ref{sec:egs:mat}) in case $\ll=2$. 

Recall that in the case of monstrous moonshine the genus zero property---that each $T_g$ should induce an isomorphism $X_{\G}\to \hat{\CC}$ (cf. (\ref{eqn:sums:XG})) for some group $\G<\SL_2(\RR)$---was the primary tool for predicting the McKay--Thompson series, and we have seen in \S\ref{sec:egs:mon} that this is equivalent to the property that $T_g$ coincide (up to an additive constant, cf. (\ref{eqn:egs:mon:TR})) with the Rademacher sum $R^{[-1]}_{\G,1,0}$ for some $\G$. In the case of Mathieu moonshine we have seen that each $H_g$ may recovered as $-2R^{[-1/8]}_{\G,\rho_g\e^{-3},1/2}$ for a suitable character $\rho_g$, and this is evidently a powerful analogue of the genus zero property of monstrous moonshine. Analogously, in the case of umbral moonshine it is conjectured \cite{UM} that each umbral McKay--Thompson series $H^{(\ll)}_g$ is recovered from a vector-valued Rademacher sum according to
\begin{gather}\label{eqn:egs:umb:HR}
	H^{(\ll)}_g=-2R^{[\m]}_{\G_0(n_g),\psi^{(\ll)}\rho^{(\ll)}_g,1/2}
\end{gather}
where $\mu=\mu^{(\ll)}$ and $\psi^{(\ll)}$ are as before, $\rho^{(\ll)}_g$ is a suitably defined (matrix-valued) function on $\G_0(n_g)$ and $n_g$ is a suitably chosen integer. (We refer to \cite[\S4.8]{UM} for more details on $\rho^{(\ll)}_g$ and $n_g$.) The conjectural identity (\ref{eqn:egs:umb:HR}) was the primary tool used in determining the concrete expressions for the $H^{(\ll)}_g$ that were furnished in \cite{UM}.

\section{Physical Applications}\label{sec:phys}

In the previous sections we have described the Rademacher summing procedure that produces a (mock) modular form by computing a certain reguralised sum over the representatives of the cosets $\G_\infty\backslash \G$, where $\G < \SL_2(\RR)$ is the modular group and $\G_\infty$ is its subgroup fixing the infinite cusp. These (mock) modular forms are often closely related to the partition function or the twisted partition function of certain two-dimensional  conformal field theories in physics. Hence, one might wonder if the associated Rademacher sum also has a physical meaning. The answer to this question is positive and in fact constituted an important part of the motivation to explore the relation between moonshine and Rademacher sums \cite{DunFre_RSMG,Cheng2011,CheDun_M24MckAutFrms}. 
 
A compelling physical interpretation of the Rademacher sum is provided by the so-called {AdS/CFT correspondence} \cite{MaldacenaAdv.Theor.Math.Phys.2:231-2521998} (also referred to as the {gauge/gravity duality} or the {holographic duality} in more general contexts), which asserts, among many other things, that the partition function of a given two dimensional CFT ``with an AdS dual" equals the partition function of another physical theory in three Euclidean dimensions with gravitational interaction and with asymptotically anti de Sitter (AdS) boundary condition. 
The correspondence, when applicable, provides both deep intuitive insights and powerful computational tools for the study of the theory.
From the fact that the only smooth three-manifold with asymptotically AdS torus boundary condition is a solid torus, it follows that the saddle points of such a partition function are labeled by the different possible ways to ``fill in the torus;" that is, the different choices of primitive cycle on the boundary torus which may become contractible in a solid torus that fills it \cite{MaldacenaJHEP9812:0051998}. These different saddle points are therefore labeled by the coset space $\G_\inf\backslash \G$, where $\G=\SL_2(\ZZ)$ \cite{Dijkgraaf2007}. From a bulk, gravitational point of view, the group $\SL_2(\ZZ)$ has an interpretation as the group of large diffeomorphisms, and $\G_\inf$ is the subgroup that leaves the contractible cycle invariant and therefore can be described by a mere change of coordinates. 
Such considerations underlie the previous use of Rademacher sums in the physics literature \cite{Dijkgraaf2007,Moore2007,BoerJHEP0611:0242006,KrausJHEP0701:0022007,Denef2007,Manschot2007,Murthy:2009dq}. 
See also \cite{Murthy:2011dk} for a refinement of this interpretation using localisation techniques. 

In the presence of a discrete symmetry of the conformal field theory theory, apart from the partition function one can also compute the twisted (or equivariant) partition function. In more details, recall that the partition function computes the dimension of the Hilbert space graded by the basic charges (the energy, for instance) of the theory. In the presence of a discrete symmetry whose action on the Hilbert space commutes with the operators associated with the basic conserved charges, more refined information can be gained by studying the twisted partition function (a trace over the Hilbert space with a group element inserted) which computes the graded group characters of the Hilbert space.
In the Lagrangian formulation of quantum field theories this twisting corresponds to a modification of the boundary condition. 
For a two dimensional CFT with an AdS gravity dual, this translates into a corresponding modification of the boundary condition in the gravitational path integral by an insertion of a group element $g$, which changes the set of  allowed saddle points. as a result, the allowed large diffeomorphisms is now given by a  discrete group $\G_g \subset \SL_2(\RR)$, generally different from $\SL_2(\ZZ)$. 

Note that when $\G \not\subset \SL_2(\ZZ)$, in particular when $\G = \G_0(n|h)+S$ where $S$ in a non-trivial subgroup of the group of exact divisors of $n/h$  (see \cite{DunFre_RSMG} for details), the above interpretation suggests that certain orbifold geometries should be included in the path integral as well as smooth geometries. We do not have a precise understanding from the gravity viewpoint as for when these extra contributions should be included. Some interpretation in terms of a $\ZZ/n\ZZ$-generalisation of the spin structure ($n=2$) have been put forward in \cite{DunFre_RSMG}. See also \cite{MalWit_QGPtnFn3D} for a related discussion. We hope further developments will shed light on this question in the future.

We have explained above how the sum over $\Gamma_\infty\backslash \G$ for $\G=\SL_2(\ZZ)$ can be thought of as a sum over the smooth, asymptotically $AdS_3$ geometries. Moreover, recent progress in the exact computation of path integrals in quantum gravity in AdS backgrounds suggests that the precise form of the regulator itself is also natural from the gravitational viewpoint. Recall the use of the Lipschitz summation formula \eq{eqn:fun:Lipsum} in reducing the Rademacher sum (\ref{eqn:sums:SL2wt2}) to a sum (\ref{eqn:sums:RK4}) of sums over (representatives of the non-trivial) double cosets of $\G_{\inf}$ in $\G$. This procedure can be applied quite generally and verifies the relationship (\ref{eqn:series:RadFou}) between Fourier coefficients of Rademacher sums and the Rademacher series. In practice then, instead of a sum over a pair of co-prime integers $(c,d)$ we can write a Rademacher sum as a generating function of sums over a single integer $c$.  This readily renders the following form for the Fourier coefficient $c_{\Gamma}(\Delta)$ of the term $q^{\Delta}$ in the Rademacher sum.
It is the infinite sum
\begin{gather}
c_{\Gamma}(\Delta) = \sum_{c=1}^\infty c_{\Gamma}(\Delta;c) ,
\end{gather}
where $ c_{\Gamma}(\Delta;c)$ takes the form of a product of a modified Bessel function with argument $\pi \sqrt{\Delta}/c$ and a Kloosterman sum (cf. \eq{eqn:series:cdef}).

In \cite{Dabholkar:2011ec} an example has been provided where the gravity path integral is argued to localise on configurations giving precisely the contribution of the above form to the gravity partition function.
First, the sum over $c$ has the interpretation as a sum over gravitational instantons obtained from orbifolding the configuration corresponding to $c=1$ by a symmetry group ${\cal G} \cong {\mathbb Z}/{c \mathbb Z}$. Second, the Bessel function arises naturally as the result of the finite-dimensional integral obtained from localising the infinite-dimensional path integral on the given instanton configuration. This result is argued to be independent of the details of the orbifold and depends only on the order $c$ of the symmetry. Finally, the Kloosterman sum and the extra numerical factor is speculated to arise from summing over different possibilities of order $c$ orbifold group ${\cal G} \cong {\mathbb Z}/{c \mathbb Z}$.
It would be very interesting to see if further developments in localising the gravity path integral will lead to a more complete understanding of quantum gravity utilising Rademacher sums.

\appendix

\section{Special Functions}\label{sec:fun}

The {\em Bernoulli numbers} $B_m$ may be defined by the following Taylor expansion.
\begin{gather}\label{eqn:fun:Ber}
	\frac{t}{e^t-1}=\sum_{m\geq 0}B_m\frac{t^m}{m!}
\end{gather}

The {\em Gamma function} $\G(s)$ and {\em lower incomplete Gamma function} $\g(s,x)$ are defined by the integrals
\begin{gather}
	\G(s)=\int_0^{\inf}t^{s-1}e^{-t}{\rm d}t\label{eqn:fun:gam}\\
	\g(s,x)=\int_0^xt^{s-1}e^{-t}{\rm d}t\label{eqn:fun:lowgam}
\end{gather}
for $s$ real and positive. The expression (\ref{eqn:fun:lowgam}) is well defined for positive real $x$ but this situation can be improved, for integration by parts yields the recurrence relation 
\begin{gather}
\g(s,x)=(s-1)\g(s-1,x)-x^{s-1}e^{-x},
\end{gather}
and this in turn leads to a power series expansion
\begin{gather}\label{eqn:fun:lowgamser}
	\g(s,x)=\frac{\G(s)}{e^x}\sum_{n\geq 0}\frac{x^{n+s}}{\G(n+s+1)}
\end{gather}
which converges absolutely and locally uniformly for $x$ in $\CC$. 

For the exponential $x^s$ we employ the principal branch of the logarithm, so that 
\begin{gather}\label{eqn:fun:pbranch}
	x^s=|x|^se^{\ii\th s}\quad{\rm whenever}\quad
	x=|x|e^{\ii\th},\,-\pi<\th\leq \pi.
\end{gather}

The {\em modified Bessel function of the first kind} is denoted $I_{\a}(x)$ and may be defined by the power series expression
\begin{gather}\label{eqn:fun:Bes}
	I_{\a}(z)=\sum_{n\geq 0}\frac{1}{\G(m+\a+1)m!}\left(\frac{z}{2}\right)^{2m+\a}
\end{gather}
which converges absolutely and locally uniformly in $z$ so long as $z$ avoids the negative reals (cf. (\ref{eqn:fun:pbranch})). We consider only non-negative real values of $\a$ in this article.

The {\em Dedekind eta function}, denoted $\eta(\t)$, is a holomorphic function on the upper half-plane defined by the infinite product 
\begin{gather}\label{eqn:fun:Ded}
	\eta(\t)=q^{1/24}\prod_{n>0}(1-q^n)
\end{gather}
where $q=\ex(\t)=e^{\tpi \t}$. It is a modular form of weight $1/2$ for the modular group $\SL_2(\ZZ)$ with multiplier $\e:\SL_2(\ZZ)\to\CC$ so that
\begin{gather}\label{eqn:fun:eps}
	\eta(\g\t)\e(\g)\jac(\g,\t)^{1/4}=\eta(\t)
\end{gather}
for all $\g = \big(\begin{smallmatrix} a&b\\ c&d \end{smallmatrix}\big) \in\SL_2(\ZZ)$, where $\jac(\g,\t)=(c\t+d)^{-2}$. The {\em multiplier system} $\e$ may be described explicitly as 
\begin{gather}\label{eqn:fun:dedmlt}
\e\bem a&b\\ c&d\eem 
	=
\begin{cases}
	\ex(-b/24),&c=0,\,d=1\\
	\ex(-(a+d)/24c+s(d,c)/2+1/8),&c>0
\end{cases}
\end{gather}
where $s(d,c)=\sum_{m=1}^{c-1}(d/c)((md/c))$ and $((x))$ is $0$ for $x\in\ZZ$ and $x-\lfloor x\rfloor-1/2$ otherwise. We can deduce the values $\e(a,b,c,d)$ for $c<0$, or for $c=0$ and $d=-1$, by observing that $\e(-\g)=\e(\g)\ex(1/4)$ for $\g\in\SL_2(\ZZ)$. Observe that
\begin{gather}\label{eqn:fun:dedmlta}
\e(T^m\g)=\e(\g T^m)=\ex(-m/24)\e(\g)
\end{gather}
for $m\in\ZZ$.

Setting $q=\ex(\t)$ and $y=\ex(z)$ we use the following conventions for the four standard {\em Jacobi theta functions}.
\begin{gather}
\begin{split}\label{eqn:fun:Jactht}
\th_1(\t,z) &= -i q^{1/8} y^{1/2} \prod_{n=1}^\inf (1-q^n) (1-y q^n) (1-y^{-1} q^{n-1})\\
\th_2(\t,z) &=  q^{1/8} y^{1/2} \prod_{n=1}^\inf (1-q^n) (1+y q^n) (1+y^{-1} q^{n-1})\\
\th_3(\t,z) &=  \prod_{n=1}^\inf (1-q^n) (1+y \,q^{n-1/2}) (1+y^{-1} q^{n-1/2})\\
\th_4(\t,z) &=  \prod_{n=1}^\inf (1-q^n) (1-y \,q^{n-1/2}) (1-y^{-1} q^{n-1/2})
\end{split}
\end{gather}

We write $\m(\t,z)$ for the {\em Appell-Lerch sum} defined by setting
\begin{gather}\label{eqn:fun:AppLer}
\m(\t,z) = \frac{-i y^{1/2}}{\th_{1}(\t,z)}\,\sum_{\ell=-\inf}^\inf \frac{(-1)^{\ell} y^n q^{\ell(\ell+1)/2}}{1-y q^\ell}.
\end{gather}

The {\em Lipschitz summation formula} is the identity 
\be\label{eqn:fun:Lipsum}
\frac{(-2\p i)^s}{\G(s)}\sum_{k=1}^\inf {(k-\alpha)^{s-1}}\ex((k-\alpha)\t)=\sum_{\ell\in\ZZ}\ex(\alpha \ell)(\t+\ell)^{-s},
\ee
valid for $\Re(s)>1$ and $0\leq \alpha <1$, where $\ex(x)=e^{\tpi x}$. A nice proof of this using Poisson summation appears in \cite{KnoRob_RieFnlEqnLipSum}. Observe that both sides of (\ref{eqn:fun:Lipsum}) converge absolutely and uniformly in $\t$ on compact subsets of $\HH$. For applications to Rademacher sums of weight less than $1$ one requires an extension of (\ref{eqn:fun:Lipsum}) to $s=1$. 
Absolute convergence on the right hand side breaks down at this point but we may consider the following useful analogue. The reader may consult \cite[\S C]{Cheng2011}, for example, for a proof of (\ref{eqn:fun:Lip_s1}), and may see \cite[\S3.3]{DunFre_RSMG} for a proof of (\ref{eqn:fun:Lip_s1a0}).
\begin{lem}\label{lem:fun:LipSumAnlg}
For $0<\alpha<1$ we have 
\be\label{eqn:fun:Lip_s1}
\sum_{k=1}^\inf\ex((k-\alpha)\t)=
\sum_{-K< \ell< K}\ex(\alpha \ell)(-2\pi i)^{-1}(\t+\ell)^{-1}
+E_K(\t)
\ee
where $E_K(\t)={\cal O}(1/K^2)$, locally uniformly for $\t\in\HH$. For $\a=0$ we have
\begin{gather}\label{eqn:fun:Lip_s1a0}
\frac{1}{2}+\sum_{k>0}\ex(k\t)=\lim_{K\to \inf}\sum_{-K<\ell<K}(-2\pi\ii)^{-1}(\t+\ell)^{-1}.
\end{gather}
\end{lem}

\addcontentsline{toc}{section}{References}
\bibliographystyle{alpha}
\bibliography{revrad_tex}

\newcommand{\etalchar}[1]{$^{#1}$}
\begin{thebibliography}{dBCD{\etalchar{+}}06}

\bibitem[And81]{And_RamLstNbk_III}
George~E. Andrews.
\newblock Ramunujan's ``lost'' notebook. {III}. {T}he {R}ogers-{R}amanujan
  continued fraction.
\newblock {\em Adv. in Math.}, 41(2):186--208, 1981.

\bibitem[Apo90]{Apo_MdlFnsDirSerNumThy}
Tom~M. Apostol.
\newblock {\em Modular functions and {D}irichlet series in number theory},
  volume~41 of {\em Graduate Texts in Mathematics}.
\newblock Springer-Verlag, New York, second edition, 1990.

\bibitem[BO06]{BringmannOno2006}
Kathrin Bringmann and Ken Ono.
\newblock The {$f(q)$} mock theta function conjecture and partition ranks.
\newblock {\em Invent. Math.}, 165(2):243--266, 2006.

\bibitem[BO07]{BriOno_ArtPrpCoeHlfIntWgtMaaPoiSrs}
Kathrin Bringmann and Ken Ono.
\newblock Arithmetic properties of coefficients of half-integral weight
  {M}aass-{P}oincar\'e series.
\newblock {\em Math. Ann.}, 337(3):591--612, 2007.

\bibitem[BO10]{BringmannOno2010}
Kathrin Bringmann and Ken Ono.
\newblock Dyson's ranks and {M}aass forms.
\newblock {\em Ann. of Math. (2)}, 171(1):419--449, 2010.

\bibitem[BO12]{BriOno_CoeffHmcMaaFrms}
Kathrin Bringmann and Ken Ono.
\newblock Coefficients of harmonic maass forms.
\newblock In Krishnaswami Alladi and Frank Garvan, editors, {\em Partitions,
  q-Series, and Modular Forms}, volume~23 of {\em Developments in Mathematics},
  pages 23--38. Springer New York, 2012.

\bibitem[Bor92]{borcherds_monstrous}
Richard~E. Borcherds.
\newblock {Monstrous moonshine and monstrous Lie superalgebras}.
\newblock {\em Invent. Math.}, 109, No.2:405--444, 1992.

\bibitem[Car07]{carnahan}
Scott Carnahan.
\newblock {\em {Monstrous Lie algebras and generalized moonshine}}.
\newblock PhD thesis, University of California, Berkeley, 2007.

\bibitem[CC11]{ChoCho_ZagDtyHmcMaaFmsIntWgt}
Bumkyu Cho and Youngju Choie.
\newblock Zagier duality for harmonic weak {M}aass forms of integral weight.
\newblock {\em Proc. Amer. Math. Soc.}, 139(3):787--797, 2011.

\bibitem[CD11]{Cheng2011}
Miranda C.~N. Cheng and John F.~R. Duncan.
\newblock {On Rademacher Sums, the Largest Mathieu Group, and the Holographic
  Modularity of Moonshine}.
\newblock October 2011.

\bibitem[CD12]{CheDun_M24MckAutFrms}
Miranda C.~N. Cheng and John F.~R. Duncan.
\newblock {The Largest Mathieu Group and (Mock) Automorphic Forms}.
\newblock January 2012.

\bibitem[CDH12]{UM}
Miranda C.~N. Cheng, John F.~R. Duncan, and Jeffrey~A. Harvey.
\newblock {Umbral Moonshine}.
\newblock 2012.

\bibitem[Che10]{Cheng2010_1}
Miranda C.~N. Cheng.
\newblock {$K3$} {S}urfaces, {$N=4$} {D}yons, and the {M}athieu {G}roup
  {$M_{24}$}.
\newblock May 2010.

\bibitem[CMS04]{ConMcKSebDiscGpsM}
John Conway, John McKay, and Abdellah Sebbar.
\newblock On the discrete groups of {M}oonshine.
\newblock {\em Proc. Amer. Math. Soc.}, 132:2233--2240, 2004.

\bibitem[CN79]{conway_norton}
J.~H. Conway and S.~P. Norton.
\newblock {Monstrous Moonshine}.
\newblock {\em Bull. London Math. Soc.}, 11:308~339, 1979.

\bibitem[Cum04]{Cum_CngSbsGpsCmmMdlGrpGns01}
C.~J. Cummins.
\newblock Congruence subgroups of groups commensurable with {${\rm PSL}(2,\Bbb
  Z)$} of genus 0 and 1.
\newblock {\em Experiment. Math.}, 13(3):361--382, 2004.

\bibitem[dBCD{\etalchar{+}}06]{BoerJHEP0611:0242006}
Jan de~Boer, Miranda C.~N. Cheng, Robbert Dijkgraaf, Jan Manschot, and Erik
  Verlinde.
\newblock A farey tail for attractor black holes.
\newblock {\em JHEP}, 0611:024,2006, JHEP0611:024,2006.

\bibitem[DF11]{DunFre_RSMG}
John F.~R. Duncan and Igor~B. Frenkel.
\newblock Rademacher sums, moonshine and gravity.
\newblock {\em Commun. Number Theory Phys.}, 5(4):1--128, 2011.

\bibitem[DGM11]{Dabholkar:2011ec}
Atish Dabholkar, Joao Gomes, and Sameer Murthy.
\newblock {Localization and Exact Holography}.
\newblock 2011.

\bibitem[DM11]{Denef2007}
Frederik Denef and Gregory~W. Moore.
\newblock {Split states, entropy enigmas, holes and halos}.
\newblock {\em JHEP}, 1111:129, 2011.
\newblock 149 pages, 21 figures.

\bibitem[DMMV07]{Dijkgraaf2007}
Robbert Dijkgraaf, Juan Maldacena, Gregory Moore, and Erik Verlinde.
\newblock A black hole farey tail.
\newblock 2007.

\bibitem[DMZ]{Atish_Sameer_Don}
Atish Dabholkar, Sameer Murthy, and Don Zagier.
\newblock {Quantum Black Holes, Wall-Crossing, and Mock Modular Forms}.
\newblock {\em to appear}.

\bibitem[EH09]{Eguchi2009a}
Tohru Eguchi and Kazuhiro Hikami.
\newblock {Superconformal Algebras and Mock Theta Functions 2. Rademacher
  Expansion for K3 Surface}.
\newblock {\em Communications in Number Theory and Physics}, 3,:531--554, April
  2009.

\bibitem[EH11]{Eguchi2010a}
Tohru Eguchi and Kazuhiro Hikami.
\newblock {Note on Twisted Elliptic Genus of K3 Surface}.
\newblock {\em Phys.Lett.}, B694:446--455, 2011.

\bibitem[EOT11]{Eguchi2010}
Tohru Eguchi, Hirosi Ooguri, and Yuji Tachikawa.
\newblock {Notes on the K3 Surface and the Mathieu group $M_{24}$}.
\newblock {\em Exper.Math.}, 20:91--96, 2011.

\bibitem[FLM84]{FLMPNAS}
I.~B. Frenkel, J.~Lepowsky, and A.~Meurman.
\newblock A natural representation of the {F}ischer-{G}riess {M}onster with the
  modular function {$J$} as character.
\newblock {\em Proc. Nat. Acad. Sci. U.S.A.}, 81(10, Phys. Sci.):3256--3260,
  1984.

\bibitem[FLM85]{FLMBerk}
Igor~B. Frenkel, James Lepowsky, and Arne Meurman.
\newblock A moonshine module for the {M}onster.
\newblock In {\em Vertex operators in mathematics and physics (Berkeley,
  Calif., 1983)}, volume~3 of {\em Math. Sci. Res. Inst. Publ.}, pages
  231--273. Springer, New York, 1985.

\bibitem[FLM88]{FLM}
Igor Frenkel, James Lepowsky, and Arne Meurman.
\newblock {\em Vertex operator algebras and the {M}onster}, volume 134 of {\em
  Pure and Applied Mathematics}.
\newblock Academic Press Inc., Boston, MA, 1988.

\bibitem[FO08]{FolOno_DltyMckThtFn}
Amanda Folsom and Ken Ono.
\newblock Duality involving the mock theta function {$f(q)$}.
\newblock {\em J. Lond. Math. Soc. (2)}, 77(2):320--334, 2008.

\bibitem[GHV10a]{Gaberdiel2010a}
Matthias~R. Gaberdiel, Stefan Hohenegger, and Roberto Volpato.
\newblock {Mathieu Moonshine in the elliptic genus of K3}.
\newblock {\em JHEP}, 1010:062, 2010.

\bibitem[GHV10b]{Gaberdiel2010}
Matthias~R. Gaberdiel, Stefan Hohenegger, and Roberto Volpato.
\newblock {Mathieu twining characters for K3}.
\newblock {\em JHEP}, 1009:058, 2010.
\newblock 19 pages.

\bibitem[GS83]{GolSar_Kloo}
D.~Goldfeld and P.~Sarnak.
\newblock Sums of {K}loosterman sums.
\newblock {\em Invent. Math.}, 71(2):243--250, 1983.

\bibitem[Gue09]{Gue_HmcMaaMdlGdsEvnIntWts}
P.~Guerzhoy.
\newblock On weak harmonic {M}aass-modular grids of even integral weights.
\newblock {\em Math. Res. Lett.}, 16(1):59--65, 2009.

\bibitem[KL07]{KrausJHEP0701:0022007}
Per Kraus and Finn Larsen.
\newblock Partition functions and elliptic genera from supergravity.
\newblock {\em JHEP}, 0701:002,2007, JHEP 0701:002,2007.

\bibitem[Kno61a]{Kno_ConstMdlrFnsI}
Marvin~Isadore Knopp.
\newblock Construction of a class of modular functions and forms.
\newblock {\em Pacific J. Math.}, 11:275--293, 1961.

\bibitem[Kno61b]{Kno_ConstMdlrFnsII}
Marvin~Isadore Knopp.
\newblock Construction of a class of modular functions and forms. {II}.
\newblock {\em Pacific J. Math.}, 11:661--678, 1961.

\bibitem[Kno62a]{Kno_ConstAutFrmsSuppSeries}
Marvin~Isadore Knopp.
\newblock Construction of automorphic forms on {$H$}-groups and supplementary
  {F}ourier series.
\newblock {\em Trans. Amer. Math. Soc.}, 103:168--188, 1962.

\bibitem[Kno62b]{Kno_AbIntsMdlrFns}
Marvin~Isadore Knopp.
\newblock On abelian integrals of the second kind and modular functions.
\newblock {\em Amer. J. Math.}, 84:615--628, 1962.

\bibitem[Kno86]{Kno_SmlPosPowTheta}
M.~I. Knopp.
\newblock On the {F}ourier coefficients of small positive powers of {$\theta
  (\tau)$}.
\newblock {\em Invent. Math.}, 85(1):165--183, 1986.

\bibitem[Kno89]{Kno_SmlPosWgt}
Marvin~I. Knopp.
\newblock On the {F}ourier coefficients of cusp forms having small positive
  weight.
\newblock In {\em Theta functions---{B}owdoin 1987, {P}art 2 ({B}runswick,
  {ME}, 1987)}, volume~49 of {\em Proc. Sympos. Pure Math.}, pages 111--127.
  Amer. Math. Soc., Providence, RI, 1989.

\bibitem[Kno90]{Kno_RadonJPoinSerNonPosWtsEichCohom}
Marvin~I. Knopp.
\newblock Rademacher on {$J(\tau),$} {P}oincar\'e series of nonpositive weights
  and the {E}ichler cohomology.
\newblock {\em Notices Amer. Math. Soc.}, 37(4):385--393, 1990.

\bibitem[Kow10]{Kow_PoiAncNmbThy}
Emmanuel Kowalski.
\newblock Poincar\'e and analytic number theory.
\newblock In {\em The scientific legacy of {P}oincar\'e}, volume~36 of {\em
  Hist. Math.}, pages 73--85. Amer. Math. Soc., Providence, RI, 2010.

\bibitem[KR01]{KnoRob_RieFnlEqnLipSum}
Marvin Knopp and Sinai Robins.
\newblock Easy proofs of {R}iemann's functional equation for {$\zeta(s)$} and
  of {L}ipschitz summation.
\newblock {\em Proc. Amer. Math. Soc.}, 129(7):1915--1922 (electronic), 2001.

\bibitem[Mal98]{MaldacenaAdv.Theor.Math.Phys.2:231-2521998}
Juan~M. Maldacena.
\newblock The large {N} limit of superconformal field theories and
  supergravity.
\newblock {\em Adv.Theor.Math.Phys.}, 2:231--252, 1998.

\bibitem[MM07]{Manschot2007}
Jan Manschot and Gregory~W. Moore.
\newblock A modern fareytail.
\newblock {\em Commun.Num.Theor.Phys.}, 4,:103--159,2010, December 2007.

\bibitem[MN11]{Murthy:2011dk}
Sameer Murthy and Satoshi Nawata.
\newblock {Which AdS3 Configurations Contribute to the SCFT2 Elliptic Genus?}
\newblock 2011.

\bibitem[Moo07]{Moore2007}
Gregory~W. Moore.
\newblock Les {H}ouches {L}ectures on {S}trings and {A}rithmetic.
\newblock 2007.

\bibitem[MP09]{Murthy:2009dq}
Sameer Murthy and Boris Pioline.
\newblock {A Farey tale for N=4 dyons}.
\newblock {\em JHEP}, 0909:022, 2009.

\bibitem[MS98]{MaldacenaJHEP9812:0051998}
Juan Maldacena and Andrew Strominger.
\newblock Ads3 black holes and a stringy exclusion principle.
\newblock {\em JHEP}, 9812:005,1998, JHEP 9812:005,1998.

\bibitem[MW10]{MalWit_QGPtnFn3D}
Alexander Maloney and Edward Witten.
\newblock Quantum gravity partition functions in three dimensions.
\newblock {\em J. High Energy Phys.}, (2):029, 58, 2010.

\bibitem[Nie74]{Nie_ConstAutInts}
Douglas Niebur.
\newblock Construction of automorphic forms and integrals.
\newblock {\em Trans. Amer. Math. Soc.}, 191:373--385, 1974.

\bibitem[Pet30]{Pet_AutFrmDtgArtPoiRhn}
Hans Petersson.
\newblock Theorie der automorphen {F}ormen beliebiger reeller {D}imension und
  ihre {D}arstellung durch eine neue {A}rt {P}oincar\'escher {R}eihen.
\newblock {\em Math. Ann.}, 103(1):369--436, 1930.

\bibitem[Pet32]{Pet_UbrEntAutFrm}
Hans Petersson.
\newblock \"{U}ber die {E}ntwicklungskoeffizienten der automorphen {F}ormen.
\newblock {\em Acta Math.}, 58(1):169--215, 1932.

\bibitem[Pet33]{Pet_UbrEntAllKlasAutFrm}
Hans Petersson.
\newblock \"{U}ber die {E}ntwicklungskoeffizienten einer allgemeinen {K}lasse
  automorpher {F}ormen.
\newblock {\em Math. Ann.}, 108(1):370--377, 1933.

\bibitem[Poi11]{Poi_FtnMdlFtnFuc}
H.~Poincar{\'e}.
\newblock Fonctions modulaires et fonctions fuchsiennes.
\newblock {\em Ann. Fac. Sci. Toulouse Sci. Math. Sci. Phys. (3)}, 3:125--149,
  1911.

\bibitem[Pri99]{Pri_SmlPosWgt_I}
Wladimir de~Azevedo Pribitkin.
\newblock The {F}ourier coefficients of modular forms and {N}iebur modular
  integrals having small positive weight. {I}.
\newblock {\em Acta Arith.}, 91(4):291--309, 1999.

\bibitem[Pri00a]{Pri_SmlPosWgt_II}
Wladimir de~Azevedo Pribitkin.
\newblock The {F}ourier coefficients of modular forms and {N}iebur modular
  integrals having small positive weight. {II}.
\newblock {\em Acta Arith.}, 93(4):343--358, 2000.

\bibitem[Pri00b]{Pri_GnlzdGolSarEst}
Wladimir de~Azevedo Pribitkin.
\newblock A generalization of the {G}oldfeld-{S}arnak estimate on {S}elberg's
  {K}loosterman zeta-function.
\newblock {\em Forum Math.}, 12(4):449--459, 2000.

\bibitem[Rad37]{Rad_PtnFn}
Hans Rademacher.
\newblock On the {P}artition {F}unction $p(n)$.
\newblock {\em Proc. London Math. Soc. (2)}, 43:241--254, 1937.

\bibitem[Rad38]{Rad_FouCoeffMdlrInv}
Hans Rademacher.
\newblock The {F}ourier {C}oefficients of the {M}odular {I}nvariant
  {J}({$\tau$}).
\newblock {\em Amer. J. Math.}, 60(2):501--512, 1938.

\bibitem[Rad39]{Rad_FuncEqnModInv}
Hans Rademacher.
\newblock The {F}ourier {S}eries and the {F}unctional {E}quation of the
  {A}bsolute {M}odular {I}nvariant {J}({$\tau$}).
\newblock {\em Amer. J. Math.}, 61(1):237--248, 1939.

\bibitem[Rog17]{Rog_FalseTheta}
L.~J. Rogers.
\newblock On two theorems of combinatory analysis and some allied identities.
\newblock {\em Proc. London Math. Soc.}, 16:315--316, 1917.

\bibitem[Rou06]{Rou_ZagDtyExpBorPdtHilMdrFrms}
Jeremy Rouse.
\newblock Zagier duality for the exponents of {B}orcherds products for
  {H}ilbert modular forms.
\newblock {\em J. London Math. Soc. (2)}, 73(2):339--354, 2006.

\bibitem[RZ38]{RadZuc_FouCoeffMdlrFrms}
Hans Rademacher and Herbert~S. Zuckerman.
\newblock On the {F}ourier coefficients of certain modular forms of positive
  dimension.
\newblock {\em Ann. of Math. (2)}, 39(2):433--462, 1938.

\bibitem[Sel65]{Sel_EstFouCoeffs}
Atle Selberg.
\newblock On the estimation of {F}ourier coefficients of modular forms.
\newblock In {\em Proc. {S}ympos. {P}ure {M}ath., {V}ol. {VIII}}, pages 1--15.
  Amer. Math. Soc., Providence, R.I., 1965.

\bibitem[Shi71]{Shi_IntThyAutFns}
Goro Shimura.
\newblock {\em Introduction to the arithmetic theory of automorphic functions}.
\newblock Publications of the Mathematical Society of Japan, No. 11. Iwanami
  Shoten, Publishers, Tokyo, 1971.
\newblock Kan{\^o} Memorial Lectures, No. 1.

\bibitem[Siv89]{Siv_ClsThyArtFns}
R.~Sivaramakrishnan.
\newblock {\em Classical theory of arithmetic functions}, volume 126 of {\em
  Monographs and Textbooks in Pure and Applied Mathematics}.
\newblock Marcel Dekker Inc., New York, 1989.

\bibitem[Tho79a]{Tho_FinGpsModFns}
J.~G. Thompson.
\newblock Finite groups and modular functions.
\newblock {\em Bull. London Math. Soc.}, 11(3):347--351, 1979.

\bibitem[Tho79b]{Tho_NmrlgyMonsEllModFn}
J.~G. Thompson.
\newblock Some numerology between the {F}ischer-{G}riess {M}onster and the
  elliptic modular function.
\newblock {\em Bull. London Math. Soc.}, 11(3):352--353, 1979.

\bibitem[Zag02]{Zag_TrcSngMdl}
Don Zagier.
\newblock Traces of singular moduli.
\newblock In {\em Motives, polylogarithms and {H}odge theory, {P}art {I}
  ({I}rvine, {CA}, 1998)}, volume~3 of {\em Int. Press Lect. Ser.}, pages
  211--244. Int. Press, Somerville, MA, 2002.

\bibitem[Zag09]{zagier_mock}
Don Zagier.
\newblock Ramanujan's mock theta functions and their applications (after
  {Z}wegers and {O}no-{B}ringmann).
\newblock {\em Ast\'erisque}, (326):Exp. No. 986, vii--viii, 143--164 (2010),
  2009.
\newblock S{\'e}minaire Bourbaki. Vol. 2007/2008.

\bibitem[Zwe02]{zwegers}
Sander Zwegers.
\newblock {\em {Mock Theta Functions}}.
\newblock PhD thesis, Utrecht University, 2002.

\end{thebibliography}
\end{document}